\documentclass[leqno,11pt]{amsart}
\usepackage{amssymb, amsmath,latexsym,amsfonts,amsbsy, amsthm,esint}
\usepackage{color}
\usepackage{graphicx}
\usepackage{tikz}
\usetikzlibrary{arrows, calc, decorations.markings,intersections, patterns,patterns.meta}
\usepackage{caption}
\usepackage{subcaption}
\usepackage{hyperref}
\usepackage{enumitem}
\usepackage{microtype}
\usepackage{booktabs}

\usepackage{pgfplots}
\pgfplotsset{compat=1.15}
\usepackage{mathrsfs}
\usepackage[foot]{amsaddr}

\makeatletter
\DeclareFontFamily{U}{tipa}{}
\DeclareFontShape{U}{tipa}{m}{n}{<->tipa10}{}

\newcommand{\arc@char}{%
  {\usefont{U}{tipa}{m}{n}\symbol{62}}%
}

\newcommand{\arc}[1]{\mathpalette\arc@arc{#1}}

\newcommand{\arc@arc}[2]{%
  \sbox0{$\m@th#1#2$}%
  \vbox{%
    \hbox{\resizebox{\wd0}{\height}{\arc@char}}%
    \nointerlineskip
    \box0
  }%
}
\makeatother

\let\pa\partial

\let\g\gamma
\let\d\delta

\let\s\sigma

\let\f\frac

\let\p\psi

\let\G\Gamma

\let\Om\Omega

\let\e\varepsilon

\newcommand{\BR}{\mathbb R}

\newcommand{\B}{\mathcal B}
\newcommand{\J}{\mathcal J}

\newcommand{\1}{\mathbf 1}

\DeclareMathOperator{\dist}{dist}

\DeclareMathOperator{\ang}{ang}

\DeclareMathOperator*{\esssup}{ess\,sup}
\DeclareMathOperator*{\essinf}{ess\,inf}

\newcommand{\ch}{\mathcal H^1}

\newcommand{\beq}{\begin{equation}}
\newcommand{\eeq}{\end{equation}}
\newcommand{\beqo}{\begin{equation*}}
\newcommand{\eeqo}{\end{equation*}}
\numberwithin{equation}{section}

\theoremstyle{plain}

\newtheorem{theorem}{Theorem}[section]

\newtheorem{lemma}[theorem]{Lemma}
\newtheorem{proposition}[theorem]{Proposition}
\newtheorem{corol}[theorem]{Corollary}

\theoremstyle{remark}

\newtheorem{rmk}{Remark}[section]

\begin{document}

\title[planar four-junction solution]{On the planar four-junction solution to a vector-valued Allen-Cahn system}

\author{Nicholas D. Alikakos\textsuperscript{1}}
\address{\textsuperscript{1}Department of Mathematics, University of Athens (EKPA), Panepistemiopolis, 15784 Athens, Greece}
\email{nalikako@math.uoa.gr}

\author{Giorgio Fusco\textsuperscript{2}}
\address{\textsuperscript{2}Dipartimento di Matematica Pura ed Applicata, Universit\`{a} degli Studi dell'Aquila, Via Vetoio, 67010 Coppito, L'Aquila, Italy}
\email{fusco@univaq.it}

\author{Zhiyuan Geng\textsuperscript{3}}
\address{\textsuperscript{3}Department of Mathematical Sciences, University of Arkansas, Fayetteville, AR 72701, USA}
\email{zgeng@uark.edu}

\begin{abstract}
    This paper concerns entire minimizing four-junction solutions to a two-dimensional Allen-Cahn system with a four-well potential. In a suitable regime of unequal surface tensions, the construction yields a family of solutions with varying opening angles of phase sectors. The proof employs two-direction slicing arguments for tight energy bounds and localization of the diffuse interface, together with a new iteration scheme that controls change of approximating sharp interfaces at successively smaller scales. 
\end{abstract}


\date{\today}

\maketitle

\section{Introduction}\label{sec:intro}

There is a very extensive literature developing and elaborating on the relationship between solutions of the scalar Allen-Cahn equation
\begin{equation}\label{intro:scalar-AC}
 \Delta u-W'(u)=0,\qquad u:\BR^n\to\BR,
\end{equation}
with the prototype $W(u)=\tfrac14(1-u^2)^2$, and minimal surfaces, one among many deep discoveries of De Giorgi in the 1970's. Expository works and lectures include those of Pacard \cite{pacard2012role}, del Pino, Kowalczyk and Wei \cite{delpino2012beyond}, Chodosh \cite{chodosh2023notes}, Savin \cite{savin2010survey}, and Farina and Valdinoci \cite{farinavaldinoci2009survey}. The foundational work of Modica-Mortola \cite{modicamortola1977} and Modica \cite{modica1987gradient} identifies the perimeter functional as the $\Gamma$-limit of
\[
J_\e(u,\Omega)=\int_\Omega\left(\frac\e2|\nabla u|^2+\frac1\e W(u)\right)\,dx.
\]
For the double-well potential above, the limiting energy is $\sigma\operatorname{Per}(E;\Omega)$, where $\sigma=\int_{-1}^1\sqrt{2W(t)}\,dt$. Under appropriate hypotheses, minimizers $u_\e$ of $J_\e$ converge along a subsequence to a two-phase map whose interface is the minimal surface.

By comparison, the study of the vector case
\begin{equation}\label{intro:vector-AC}
 \Delta u-W_u(u)=0,\qquad u:\BR^n\to\BR^m,\quad m\geq2,
\end{equation}
where 
$$
W\geq0,\qquad \{W=0\}=\{a_1,\ldots,a_N\},\quad a_i \text{ the phases},
$$
which is the object of the present paper, is very underdeveloped.  The basic reason behind this is the lack of the maximum principle for \eqref{intro:vector-AC} and its main consequence, the lack of the monotonicity formula, a PDE fact established in the scalar case by Modica \cite{modica1985gradient} (although there has been very recent progress
on this last issue by Bethuel \cite{bethuel2025asymptotics} and Laux and Takasao \cite{laux2026brakke}). 

Fonseca and Tartar \cite{fonseca1989gradient}, and Sternberg \cite{sternberg1988effect} in the late 1980s, and a little later Baldo \cite{baldo}, established in increasing degree of generality, the vector-valued sharp-interface limit
\begin{equation}\label{intro:general-perimeter}
 J_\e(\cdot,\Omega)\ \xrightarrow{\ \Gamma\ }\ \mathcal P_N(E;\Omega)
 :=\sum_{1\leq i<j\leq N}\sigma_{ij}\mathcal H^{n-1}
       (\partial^*E_i\cap\partial^*E_j\cap\Omega),
\end{equation}
where $E=\{E_1,\ldots,E_N\}$ is a Caccioppoli partition of $\Omega$ and the corresponding limiting BV map is
\[
 u_0=\sum_{i=1}^N a_i\mathbf1_{E_i}.
\]
Thus under suitable constraints, the minimizers $u_\e$ of $J_\e(u)$ converge, along a subsequence, to $u_0$, the map associated with a weighted minimal partition. In the Dirichlet setting considered here, the boundary values $g_\e$ on $\pa\Omega$ are prepared to take values in $\{a_1,...,a_N\}$ away from proper $O(\e)$ transition arcs. The relevant sharp-interface results are recalled in Section \ref{sec: preliminary}, with references including \cite{Gazoulis}.

In both the scalar and the vector case, the ODE solutions (heteroclinics)
\[
 \ddot U=W'(U),\qquad U(-\infty)=-1,\quad U(+\infty)=1,\quad U:\BR\to\BR
\]
and
\begin{equation}\label{intro:heteroclinic}
 \begin{gathered}
 \ddot U_{ij}=W_u(U_{ij}),\qquad U_{ij}(-\infty)=a_i,\quad U_{ij}(+\infty)=a_j,\ \ i\neq j,\\
U_{ij}\in W^{1,2}_{loc}(\BR;\BR^m),\ \ \ \sigma_{ij}=\min \int_{\BR}\left(\frac12|U_{ij}'|^2+W(U_{ij})\right)\,ds,
 \end{gathered}
\end{equation}
are fundamental. 

However, in the scalar case they are even more so, due to the Bernstein analogy that motivated De Giorgi's conjecture: at dimensions $n\leq 8$, if a solution of \eqref{intro:scalar-AC} is monotone in one direction, then it is necessarily of the form\par
\begin{equation}\label{eq degiorgi conj}
u(x)=U(d(x,\Pi))
\end{equation}
where $d$ is the signed distance function, $\Pi$ a hyperplane, and $U$ the heteroclinic solution. De Giorgi's conjecture was resolved in the celebrated works of Ghoussoub and Gui \cite{ghoussoubgui1998} for $n=2$, Ambrosio and Cabr\'e \cite{ambrosiocabre2000} for $n=3$, and Savin \cite{savin2009flat} for $4\leq n\leq8$ under an additional asymptotic hypothesis. Very recently, Chan, Fern\'andez-Real, Figalli, Florit-Simon and Serra \cite{chan2026stable} and Liu, Luo, Wang, Wei, Wei and Wu \cite{liu2026stable} have settled $n=4$ without any asymptotic condition independently. Del Pino, Kowalczyk and Wei \cite{delpino2011counterexample} constructed counterexamples for $n\geq 9$. 

On the other hand, in the vector case the corresponding objects, besides planes, are singular minimal cones. For example, for $n=2$ and $\sigma_{ij}=1$, the minimal cones are the straight line and the triod, which define the interface of the two possible optimal partitions of the plane. This goes hand in hand with the physics of three or more phases with the new phenomenon being the appearance of junctions, which are sets of higher codimension where three or more phases meet.

Proceeding in the vein of De Giorgi's analogy, one is inclined to conjecture a relationship between minimal cones and certain entire solutions of \eqref{intro:vector-AC}, which now cannot be expressed in terms of the heteroclinics since the transition at the junction is not one-dimensional, a hard problem to begin with. In light of Alberti, Ambrosio and Cabr\'e \cite{albertiambrosiocabre2001}, minimizers under compact perturbations appear to be the appropriate class.

Diffuse counterparts of the planar triod and the tetrahedral cone in $\BR^3$ were obtained first in the symmetry class of the equilateral triangle and of the tetrahedron by Bronsard, Gui and Schatzman \cite{bronsard1996three} and Gui and Schatzman \cite{gui2008quadruple} respectively, and later in Fusco \cite{fusco2024relaxed}, where the reflections in \cite{bronsard1996three} were eliminated. A general theory of equivariant solutions associated with finite reflection groups and certain lattices is developed in the book \cite{afs-book} with P. Smyrnelis and two of the authors of this paper. Finally, the entire minimizing triple junction solution without symmetry constraints was obtained in \cite{alikakos2024triple} and, independently by a distinct method, by Sandier and Sternberg \cite{ss2024}. The strategy in \cite{alikakos2024triple} was based on the method of a tight lower bound introduced in \cite{AF}.

There are also rigidity results in the direction of a vector De Giorgi conjecture for minimizers. As a byproduct in \cite{ss2024}, the first vectorial analogue for the two-dimensional minimizing Allen-Cahn solutions of two phases was obtained. It was shown that if the blow-down limit of $u$ consists of only two phases $a_i,a_j$ separated by a straight line, then $u$ must be invariant along the direction of the line. Concerning the triple junction, in \cite{geng2025rigidity} it was shown that the solution constructed possesses the stratification properties of the triod, thus confirming a De Giorgi-type conjecture for solutions possessing three phases.

We have alluded in this introduction to the fact that there is a big difference between minimal surfaces or soap bubbles, where the surface tensions are all equal, and general multiphase transitions, where the surface tension is nonconstant. Consider a symmetric surface tension matrix $\Sigma_N=(\sigma_{ij})_{i,j=1}^N$ satisfying
\begin{equation}\label{intro:tension-matrix}
 \sigma_{ii}=0,\qquad \sigma_{ij}>0\quad(i\neq j),\qquad
 \sigma_{ij}<\sigma_{ik}+\sigma_{kj}\quad(i,j,k\text{ distinct}).
\end{equation}
A cone $C_N\subset\mathbb R^2$ is a set invariant under dilations, $\lambda C_N=C_N$, $\lambda>0$, which partitions $\BR^2$ into $N$ sectors $D_1,\ldots,D_N$, with opening angles \
\[
0<\theta_i<\pi,\qquad\sum_{i=1}^N\theta_i=2\pi.
\] 
If $\tau_i$ is the outward-pointing unit vector defining the ray separating $D_i$ and $D_{i+1}$, stationarity at the junction implies
\begin{equation}\label{intro:force-balance}
 \sum_{i=1}^N\sigma_{i,i+1}\tau_i=0,
\end{equation}
with indices taken cyclically. The paired-calibration method of Lawlor and Morgan \cite{lawlor1994paired} gives a useful criterion for minimality: the cone $C_N$ is minimizing if it can be calibrated, that is if there exist vectors $\xi_1,\ldots,\xi_N\in\BR^2$ such that
\begin{equation}\label{intro:paired-calibration}
 \begin{cases}
 \xi_i-\xi_j=\sigma_{ij}n_{ij},&\text{on each interface separating sectors }i,j,\\
 |\xi_i-\xi_j|\leq\sigma_{ij},&\text{for every }i\ne j,
 \end{cases}
\end{equation}
where $n_{ij}$ is the unit normal pointing from sector $i$ into sector $j$.

These vectors form a convex polygon, the dual polygon to the cone $C_N$. The equilibrium condition \eqref{intro:force-balance} determines the closure of the polygon, while sector convexity gives its strict convexity, and \eqref{intro:paired-calibration} gives minimality. Clearly a cone is compatible with a family of surface-tension matrices. Thus phase transition networks have a much richer geometry, and can meet at junctions of arbitrary higher order, with the angles between the segments determined by the relative costs.

In the present paper we consider four phases with neighboring and opposite phases correspondingly having equal surface tension,\par
\begin{equation*}
\sigma_{12}=\sigma_{23}=\sigma_{34}=\sigma_{41}=1,\qquad \sigma_{13}=\sigma_{24}=s
\end{equation*}
thus\par
\[
\Sigma_4(s)=\begin{pmatrix}0&1&s&1\\1&0&1&s\\s&1&0&1\\1&s&1&0\end{pmatrix}.
\]
The associated cone is centrally symmetric due to the equilibrium condition \eqref{intro:force-balance}. Accordingly, we consider a \emph{cross} $\Gamma=L_1\cup L_2$, the union of two straight lines intersecting at the point $O_\Gamma$. After a rotation, the lines can be written as
\begin{equation}\label{cross lines}
 L_1=O_\Gamma+\BR(\cos\phi,-\sin\phi),\qquad
 L_2=O_\Gamma+\BR(\cos\phi,\sin\phi),\qquad
 0<\phi<\frac\pi2.
\end{equation}
The cross $\G$ divides $\BR^2$ into four sectors $D_1,D_2,D_3,D_4$, labeled counterclockwise so that $D_1$ is the sector opening to the right. Their opening angles are $2\phi,\pi-2\phi,2\phi,\pi-2\phi$, respectively. It is easy to check that $\G$ can be calibrated if and only if
\begin{equation}\label{intro:weak-cross-condition}
 \max\{\sin\phi,\cos\phi\}\leq\frac s2.
\end{equation}
It follows that $\sqrt2\leq s<2$. Our construction uses the strict version of \eqref{intro:weak-cross-condition}. The equality case however is special and is discussed in the Remark following the main theorem.

\begin{figure}[ht]
\centering
\begin{tikzpicture}[scale=2.25,line cap=round]
 \draw[gray!65] (0,0) circle (1);
 \draw[thick,blue!65!black] (-38:1)--(142:1);
 \draw[thick,blue!65!black] (38:1)--(218:1);
 \draw[dotted] (0,0)--(0.57,0);
 \draw (0.36,0) arc[start angle=0,end angle=38,radius=0.36];
 \node at (19:0.49) {$\phi$};
 \fill (0,0) circle (0.022) node[below right] {$O_\Gamma$};
 \node at (0.78,0) {$D_1$};
 \node at (0,0.74) {$D_2$};
 \node at (-0.75,0) {$D_3$};
 \node at (0,-0.74) {$D_4$};
\end{tikzpicture}
\caption{The cross $\Gamma$, its center $O_\Gamma$, and the partition $\{D_i\}_{i=1}^4$.}
\label{fig:cross}
\end{figure}
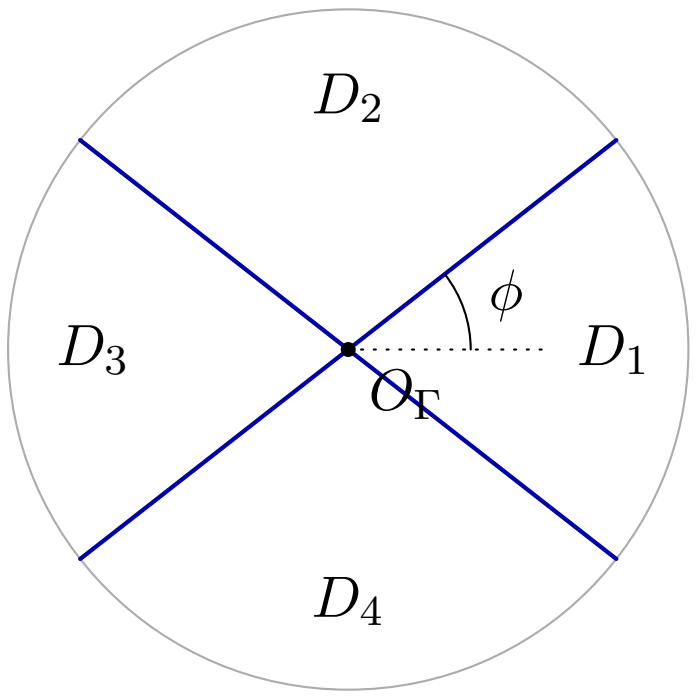

Before stating our main result we introduce the setting and the hypotheses.  For $u:\Omega\subset\BR^2\to\BR^2$, let
\begin{equation}\label{AC energy}
J_\e(u,\Omega)=\int_\Omega\left(\frac\e2|\nabla u|^2+\frac1\e W(u)\right)\,dx,
\end{equation}
with the associated Euler--Lagrange equation
\begin{equation}\label{EL eq}
 \e\Delta u-\frac1\e W_u(u)=0.
\end{equation}
We assume the following on the potential $W$.
\begin{enumerate}[label=\textup{(W\arabic*)},leftmargin=2.8em]
\itemsep0.5em
\item\label{W1}
 $W\in C^2(\BR^2;[0,\infty))$, $\{W=0\}=\{a_1,a_2,a_3,a_4\}$, and $D^2W(a_i)$ is positive definite for every $i$. Furthermore, there is $M>0$ such that
 \begin{equation}\label{radial coercivity}
 v\cdot DW(v)>0\qquad\text{whenever }|v|\geq M.
 \end{equation}
\item\label{W2}
 For the metric $d_W$ defined in \eqref{def:d_W}, set $\sigma_{ij}=d_W(a_i,a_j)$. Assume
 \begin{equation}\label{tensions}
 \sigma_{i,i+1}=1,\qquad \sigma_{13}=\sigma_{24}=s,
 \qquad s\in(\sqrt2,2).
 \end{equation}
 For every $i\neq j$, the distance $\sigma_{ij}$ is realized by a minimizing heteroclinic connection $U_{ij}\in H^1_{\mathrm{loc}}(\BR;\BR^2)\cap C^\infty(\BR;\BR^2)$:
 \[
\sigma_{ij}=\int_{\BR}\left(\frac12|U_{ij}'|^2+W(U_{ij})\right)\,dt, \qquad U_{ij}(-\infty)=a_i,\quad U_{ij}(+\infty)=a_j.
 \]
\end{enumerate}
Condition \eqref{radial coercivity} gives an $L^\infty$ bound for Dirichlet minimizers with bounded boundary values, which further implies the $C^{2,\alpha}_{\mathrm{loc}}$ regularity for any $\alpha\in(0,1)$ of $u_\varepsilon$ by standard elliptic regularity theory. For more properties of the heteroclinic $U_{ij}$ see Lemma \ref{lem:heteroclinic}.

The nondegeneracy of the wells yields a small $\delta_W>0$ such that
\begin{gather}
 4\delta_W<\min_{i\ne j}|a_i-a_j|,\label{well separation}\\
 c|p-a_i|^2\leq W(p)\leq C|p-a_i|^2
 \qquad\text{for }p\in B_{2\delta_W}(a_i).\label{near well}
\end{gather}
See also \cite{AF}. In the four-phase setting, write the sharp-interface functional as
\begin{equation}\label{def:weighted perimeter}
 \mathcal P(E;\Omega):=\sum_{1\leq i<j\leq4}\sigma_{ij}
 \mathcal H^1(\partial^*E_i\cap\partial^*E_j\cap\Omega).
\end{equation}
The partition associated with $\Gamma$ minimizes this functional under compactly supported perturbations whenever \eqref{intro:weak-cross-condition} holds, see Proposition \ref{prop:cross-calibration}. 

\begin{theorem}[Entire four-junction minimizing solution]\label{thm:main}
Assume \textup{(W1)} and \textup{(W2)}, and let $\Gamma$ be a cross centered at the origin, with half-opening-angle $\phi$ satisfying
\begin{equation}\label{strict-cross-condition}
 \max\{\sin\phi,\cos\phi\}<\frac s2.
\end{equation}
Then there exists $U\in C^2(\BR^2;\BR^2)$ solving $\Delta U=W_u(U)$ that minimizes $J_1$ in the sense of De Giorgi, that is for every open bounded set $\Omega\Subset\BR^2$,
\[
 J_1(U;\Omega)\leq J_1(V;\Omega)
 \qquad\text{whenever }V-U\in H^1_0(\Omega;\BR^2).
\]
Moreover, the full blow-down of $U$ at infinity is given by the four-partition determined by $\Gamma$:
\begin{equation}\label{eq:blowdown}
 U(Rx)\longrightarrow u_\Gamma(x):=\sum_{i=1}^4a_i\mathbf1_{D_i}(x)
 \quad\text{in }L^1_{\mathrm{loc}}(\BR^2;\BR^2)
 \quad\text{as }R\to\infty.
\end{equation}
In particular, the limiting map is independent of the sequence $R\to\infty$.
\end{theorem}
The proof also gives a quantitative rate of convergence. For a sufficiently small fixed $\e_*>0$, depending only on the potential and the angle gap $\left(\arcsin(\f{s}{2})-\max\{\phi,\f{\pi}{2}-\phi\}\right)$:
\begin{equation}\label{intro:quantitative-blowdown}
 \int_{B_1}|U(Rx)-u_\Gamma(x)|\,dx
 \leq CR^{-1/12}+C\e_*^{-2}R^{-2},\qquad R\geq4/\e_*;
\end{equation}
see \eqref{eq:large-ball-L1}.

\begin{rmk}[On the inequality \eqref{strict-cross-condition}]\label{intro:rmk-marginal}
The strict inequality in \eqref{strict-cross-condition} excludes a sharp-interface degeneracy. To explain it, consider the choice
\begin{equation}\label{intro:marginal-angle}
 \cos\phi=\frac s2, \quad \sqrt{2}\leq s<2.
\end{equation}
With the sector labeling above, the calibration vectors in \eqref{def:paired cal vec} form a rhombus with unit side lengths and diagonals
\[
 |\xi_1-\xi_3|=2\sin\phi,\qquad
 |\xi_2-\xi_4|=2\cos\phi.
\]
Thus the inequality in the second line of \eqref{intro:paired-calibration} is saturated for the opposite phases $a_2,a_4$. The diagonal $\xi_2\xi_4$ divides the rhombus into the triangles $(\xi_1,\xi_2,\xi_4)$ and $(\xi_2,\xi_3,\xi_4)$, each with side lengths $1,1,s$.

For any $d>0$, split the four-junction into two triple-junctions at $A(-d/2,0)$ and $B(d/2,0)$, joined by the horizontal segment of length $d$ separating phases $a_2$ and $a_4$. Write $\mathcal P(T_d;B_1)$ for the weighted perimeter of the partition defined by $T_d$ on $B_1$. For sufficiently small $d>0$, direct calculation gives
\begin{align*}
 \mathcal P(T_d;B_1)
 &=4\sqrt{1-\frac{d^2}{4}\sin^2\phi}+d(s-2\cos\phi)\\
 &=4-\frac12\sin^2\phi\,d^2+O(d^4),
\end{align*}
where the second equality uses \eqref{intro:marginal-angle}. Although this is smaller than $\mathcal P(\Gamma;B_1)=4$, this is not a competition partition of $\G$ due to different traces on $\pa B_1$, thus we may not directly compare their energies. It is not hard to see that $T_d$ can be calibrated by the same vectors $\xi_1,...,\xi_4$, hence it is minimizing with respect to its own trace by \cite{lawlor1994paired}. 

On the other hand we retain the segment $AB$, and connect $A$ and $B$ by straight segments to the original four boundary points of the cross $\G$ ($A$ is connected to two left points, while $B$ is connected to the right points). The resulting network $\widetilde T_d$ has the same trace as $\G$. Its energy is
\begin{align*}
 \mathcal P(\widetilde T_d;B_1)
 &=4\sqrt{1-d\cos\phi+\frac{d^2}{4}}+sd\\
 &=4+d(s-2\cos\phi)+\frac12\sin^2\phi\,d^2+O(d^3)\\
 &=4+\frac12\sin^2\phi\,d^2+O(d^3).
\end{align*}
So the cross is still minimizing under \eqref{intro:marginal-angle}, but splitting the junction only introduces quadratic cost instead of linear cost, making the four-junction degenerate. Note that the networks $T_d$ have the same four-sector blow-down as $\Gamma$. This geometry raises the possibility of a family of diffuse entire solutions with separated triple-junction cores, but neither their existence nor the existence of a minimizing four-phase solution under \eqref{intro:marginal-angle} is established here. 
\end{rmk}

These considerations above appear to be a special case of a rather general theory. At $s=\sqrt2$ and $\phi=\pi/4$, the matrix $\Sigma_4(\sqrt2)$ can be generated by taking the Euclidean distances of all pairs of vertices of a square. In general consider a regular polygon with $N$ vertices $\xi_1,...,\xi_N$ and define the surface tension matrix $\Sigma_N$ via the Euclidean distances as above. All pairwise calibration inequalities in \eqref{intro:paired-calibration} are then saturated. The dual cone associated with any subset of $k\geq3$ vertices is calibrated and remains minimizing even against competitors involving all $N$ phases, since the unused calibration vectors still satisfy the required inequalities. The various triangulations of the $N$-gon generate dual triple-junction networks, all calibrated with respect to points from the vertices set, hence all minimizing likewise. On the other hand we can introduce parameters like $s$ to cure all the degeneracies for the whole hierarchy of minimal cones, which should shadow a hierarchy of entire minimizing solutions for the vector Allen-Cahn with an $N$-well potential $W$. This is a direction for future work.

\medskip
\noindent\textbf{Main ideas of the proof.}
We begin with a minimizer $u_\e$ of $J_\e$ on the unit disk with a four-cross boundary data $g_\e$. Let $P_i$ be the four intersections of $\Gamma$ with $\partial B_1$ and $l_i=|P_i-O_\Gamma|$. First the tight lower-upper bound estimate
\begin{equation}\label{intro:energy-bounds}
 \sum_{i=1}^4l_i-C\e^{1/3}
 \leq J_\e(u_\e;B_1)
 \leq\sum_{i=1}^4l_i+C\e
\end{equation}
is derived. The upper bound follows from a direct comparison with an energy competitor, while the lower bound is obtained through a two-direction slicing argument. The horizontal and vertical energy contributions are estimated by counting the phase transitions caused by mismatch of boundary data and the interior phase configuration. Next the localization estimate 
\begin{equation}\label{intro:localization}
 |u_\e(x)-a_i|\leq C\exp\left(
 -c\frac{\max\{\dist(x,\partial(D_i\cap B_1))-C\e^{1/12},0\}}{\e}\right),
 \qquad x\in D_i\cap B_1,
\end{equation}
is established. The key ingredients here are the balanced relation between vertical and horizontal energies (Lemma \ref{lem: balance h v}) and the variational maximum principle (Proposition \ref{prop:maxp}).  Estimate \eqref{intro:localization} tells us that at a distance $O(\e^{1/12})$ from the geometric junction we can see all four diffuse phases. At this stage the localization scale $O(\e^{1/12})$ is too large to capture all four phases by a single blow-up at scale $\e$. We improve the $\e^{1/12}$ scale to $\e$ by iteration.

We work in the original variables. At step $j$, let $B_{r_j}(O_j)$ be the current disk and let $\Gamma_j$ be its associated cross, centered at $O_{\Gamma_j}=O_{j+1}$. Set
\[
 \zeta=\frac{x-O_j}{r_j},\qquad v_j(\zeta)=u_\e(O_j+r_j\zeta),\qquad\e_j=\frac{\e}{r_j},
\]
where
\[
 r_0=1,\qquad O_0=O_1=0,\qquad \Gamma_0=\Gamma,\qquad r_{j+1}=\rho_jr_j,\quad \rho_j\in[1/3,1/2].
\]
The energy scales by
\[
 J_\e(u_\e;B_{r_j}(O_j))=r_jJ_{\e_j}(v_j;B_1),
\]
so $v_j$ is still a minimizer on $B_1$ with its induced boundary values. Applying the localization estimate and transferring to the original variables yields the physical width
\[
 r_j\e_j^{1/12}=r_j\left(\frac{\e}{r_j}\right)^{1/12}=\e^{1/12}r_j^{11/12}.
\]
Thus the localization improves as the radius decreases.

Of course for this to work one needs to control the boundary conditions at each step, which clearly are not reproduced exactly. Here the localization estimate plays once more an important role. Together with an annular energy estimate, it allows us to select a smaller circle on which the boundary values are again a good trace in the sense of Section \ref{sec:good radius}. The associated sharp cross may move, thus both the centers and the orientation of interfaces will change. It turns out that the deviations are controllable:
\[
 \begin{aligned}
 |O_{\Gamma_{j+1}}-O_{\Gamma_j}|&\leq C\e^{1/12}r_j^{11/12},\\
 d_{angle}(\Gamma_{j+1},\Gamma_j)&\leq C\left(\frac{\e}{r_j}\right)^{1/12}.
 \end{aligned}
\]
Here $d_{angle}$ measures the change of the two line directions. Since the radii decrease geometrically, the total displacement
of the cross center is bounded by
\[
 |O_{\Gamma_N}-O_\Gamma| \leq C\e^{1/12}\sum_{j=0}^{N-1}r_j^{11/12} \leq C'\e^{1/12}.
\]
The angular bound might increase, but its accumulated effect satisfies
\[
 d_{angle}(\Gamma_j,\Gamma)\leq C\left(\frac{\e}{r_j}\right)^{1/12}.
\]

We stop the iteration when $r_N$ is comparable to $\e/\e_*$, where $\e_*>0$ is fixed and sufficiently small. The accumulated angular error is then at most $C\e_*^{1/12}$. We choose $\e_*$ small relative to the initial angular gap in order to keep every cross minimizing.  Meanwhile, the physical localization width at the terminal scale is $O(\e)$. We can now blow up at the final cross center $z_\e=O_{\Gamma_N}$ and let $U_\e(x)=u_\e(z_\e+\e x)$ sequentially converge to an entire minimizer $U$. The localization and angular control at all intermediate scale imply the prescribed blow-down profile for $U$ at $\infty$.

The tight lower-upper bound and the localization estimate are by now well-established ingredients in the solution of these problems, although quite challenging each time. On the other hand the iteration scheme is new and quite distinct from the backward induction scheme in \cite{alikakos2024triple}, and renders the full blow-down as $R\to\infty$, and not just a subsequence. Moreover, it is promising to apply this iteration scheme to construct entire solutions that serve as diffuse analogues of more general planar $N$-junctions and even three-dimensional junctions. 

In a forthcoming companion paper, we focus on Dirichlet problems on the disk in $\mathbb R^2$ for diffused analogs of Steiner networks, including the cross for $s=\sqrt2$, $\phi=\pi/4$.

The paper is organized as follows. Section \ref{sec: preliminary} recalls the one-dimensional estimates, the sharp-interface limit, the calibration of the cross, and the variational maximum principle. The energy bounds and the localization result are proved in Section \ref{sec:min unit disk}. Then we establish the corresponding estimates for good traces and the good-radius selection in Section \ref{sec:good radius}.  Finally, Section \ref{sec:min thm proof} carries out the iteration and proves Theorem \ref{thm:main}.

\section{Preliminaries}\label{sec: preliminary}

Throughout the rest of the paper, unless a dependence is displayed, $c,C>0$ depend only on $W$ and $\phi$. Their values may change from line to line. Let $z=(x,y)$ denote a point on the domain plane. Moreover, $|z_1-z_2|$ is the Euclidean distance, and $d_W(p,q)$ defined below is the metric in the target induced by the Allen-Cahn functional.

\subsection{The metric $d_W$ and the one-dimensional estimates}
For $p,\,q\in\mathbb{R}^2$, define
\begin{equation}\label{def:d_W}
 d_W(p,q):=\inf_{\substack{\g\in AC([0,1];\BR^2)\\ \g(0)=p,\ \g(1)=q}}\int_0^1\sqrt{2W(\g(t))}\,|\g'(t)|\,dt .
\end{equation}

For fixed $p$, the function $q\mapsto d_W(p,q)$ is locally Lipschitz and satisfies
\begin{equation}\label{metric grad}
 |D_qd_W(p,q)|=\sqrt{2W(q)} \qquad\text{for a.e. }q.
\end{equation}
\eqref{near well} implies that 
\begin{equation}\label{dW near well}
c_W|p-a_i|^2\leq d_W(a_i,p)\leq C_W|p-a_i|^2 \qquad\text{if }|p-a_i|\leq 2\delta_W.
\end{equation}

In the following lemma we list some fundamental properties of the heteroclinic connection $U_{ij}$ without proof. Readers can refer to \cite{afs-book} for detailed proofs of these properties.

\begin{lemma}\label{lem:heteroclinic}
    Let $U_{ij}$ be a heteroclinic connection between $a_i$ and $a_j$, i.e.
    \begin{equation*}
        U_{ij}=\mathrm{argmin} \int_{-\infty}^{\infty} \left(\f12|U'|^2+W(U)\right)\,dt,\quad U(-\infty)=a_i,\ \ U(+\infty)=a_j.
    \end{equation*}
    Then for every $t\in \mathbb{R}$,
    \begin{equation*}
        \f12|U'_{ij}(t)|^2=W(U_{ij}(t)).
    \end{equation*}
    Moreover, there exist constants $c,C>0$ such that 
    \begin{equation*}
    \begin{split}
        &|U_{ij}(t)-a_j|+|U_{ij}'(t)|\leq Ce^{-ct}, \quad \text{when } t\geq 0,\\
        &|U_{ij}(t)-a_i|+|U_{ij}'(t)|\leq Ce^{ct},\quad \text{when } t\leq 0.
        \end{split}
    \end{equation*}
\end{lemma}

Our next lemma provides one-dimensional energy estimates which will be heavily used in our slicing approach. 
\begin{lemma}[One-dimensional energy estimate]\label{lem:1d est}
Let \(I=(s_1,s_2)\) be an open interval and \(v\in H^1(I;\BR^2)\).  For every \(0\leq\theta\leq1\),
\begin{equation}\label{est:1d ene}
 \int_I\left(\frac{\e}{2}|v'|^2            +\frac{\theta^2}{\e}W(v)\right)\,dt\geq \theta d_W(v(s_1),v(s_2)).
\end{equation}
If \(|p-a_i|\leq \eta\) and \(|q-a_j|\leq\eta\) for some \(0<\eta<\delta_W\), then
\begin{equation}\label{est:endpoint close to ai}
 d_W(p,q)\geq\sigma_{ij}-C\eta^2.
\end{equation}
\end{lemma}

\begin{proof}
The AM-GM inequality implies
\[
 \frac{\e}{2}|v'|^2+\frac{\theta^2}{\e}W(v)
 \geq\theta\sqrt{2W(v)}\,|v'|.
\]
Integrating proves \eqref{est:1d ene}. \eqref{est:endpoint close to ai} follows directly from \eqref{dW near well} and the triangle inequality for \(d_W\).
\end{proof}

\subsection{Sharp interface limit.} 

The $\Gamma$-convergence theory has been applied to the study of the sharp interface limit of $J_\e(\cdot,\Omega)$ in \cite{baldo, fonseca1989gradient, sternberg1988effect, Gazoulis}. We summarize the results below for our need: If $u_\e\to u$ in $L^1(\Omega;\BR^2)$, then
\begin{equation}\label{G-conv:liminf}
\liminf_{\e\to0}J_\e(u_\e;\Om)
 \geq \mathcal P(\{E_i\}_{i=1}^4;\Om)
\end{equation}
whenever $u=\sum_i a_i\mathbf 1_{E_i}$. On the other hand, every finite weighted perimeter partition $\Om=\cup_{i=1}^4 E_i$ admits a recovery sequence $u_{\e_k} \to \sum_{i=1}^4 a_i\mathbf{1}_{E_i} $ in $L^1(\Om;\BR^2)$ for which equality holds as $\e_k\to 0$. Therefore, the minimal partition problem \eqref{def:weighted perimeter} is the $\Gamma$-limit associated with the Allen-Cahn functional \eqref{AC energy} as $\e\to 0$.  

For the functional \eqref{def:weighted perimeter}, the minimality of the cross $\G$ satisfying \eqref{intro:weak-cross-condition} follows by a paired calibration argument by Lawlor and Morgan \cite{lawlor1994paired}.  For completeness, we present a short proof below. Set the vectors
\begin{equation}\label{def:paired cal vec}
 \xi_1=(0,0),\qquad \xi_2=(\sin\phi,-\cos\phi),\qquad\xi_3=(2\sin\phi,0),\qquad\xi_4=(\sin\phi,\cos\phi).
\end{equation}
They satisfy
\begin{equation*}
  |\xi_i-\xi_j|\le\sigma_{ij},\qquad \forall\  i\neq j
\end{equation*}
when \eqref{intro:weak-cross-condition} holds.

\begin{proposition}[Minimality of the cross $\Gamma$] \label{prop:cross-calibration}
Let $\G$ satisfy \eqref{intro:weak-cross-condition}.  For every bounded set $\Omega$, the partition $\{D_i\}_{i=1}^4$ minimizes $\mathcal P(\cdot\,;\Omega)$ among partitions that coincide outside a compactly supported subset of $\Om$.  Equivalently, $\G$ is minimizing under every compactly supported perturbation.
\end{proposition}

\begin{proof}
Let $\{E_i\}_{i=1}^4$ be a competitor such that $E_i$ coincide with $D_i$ outside a compactly supported subset of $\Om$ for every $i$. Denote by $\nu_{ij}$ the measure-theoretical unit normal on $\pa E_i\cap \pa E_j\cap \Om$ pointing from $E_i$ to $E_j$.  Since $|\xi_i-\xi_j|\le\sigma_{ij}$, we have 
\begin{align*}
 \mathcal P(\{E_i\}_{i=1}^4;\Omega)&\geq\sum_{i<j}\int_{\partial^*E_i\cap\partial^*E_j\cap\Omega}(\xi_i-\xi_j)\cdot\nu_{ij}\,\,d\ch\\
&=\sum_i\int_{\partial^*E_i\cap\Omega}\xi_i\cdot\nu_{E_i}\,d\ch\\
&=-\sum_i \int_{E_i\cap\pa \Om} \xi_i \cdot \nu_{\overline{E_i}\cap \pa \Om} \,d\ch,
\end{align*}
The last quantity only depends on the $E_i\cap \pa\Om$, which is independent of the choice of competitors on $\Om$.  For the cross $\G$, $\xi_i-\xi_{i+1}$ is exactly the unit normal on $\pa D_i\cap\pa D_{i+1}\cap\Om$ from $D_i$ to $D_{i+1}$, thus the first inequality becomes equality for $\{D_i\}_{i=1}^4$, which proves the minimality of $\G$.
\end{proof}

\subsection{Variational maximum principle. }

We recall the following vector-valued variational maximum principle result established by Alikakos and Fusco in \cite{AF2}.
\begin{proposition}\label{prop:maxp}
Let $u$ be a minimizer of $J_\e$ in a bounded domain $\Omega$ and suppose 
\[
  |u-a_i|\leq \g\quad\text{on }\partial \Om,
  \ \text{ for some }0<\g\leq \d_W,
\]
then
\begin{equation}\label{eq:maxp}
  |u-a_i|\leq \g\quad\text{in }\Om.
\end{equation}
Moreover, there exist positive constants $c,C$, depending on $W$, such that 
\begin{equation}\label{eq:exp decay}
    |u(z)-a_i|\leq C\exp\left( -\f{c\dist(z,\pa \Om)}{\e} \right),\quad z\in\Om.
\end{equation}
\end{proposition}

\begin{proof}
The first conclusion is the variational maximum principle of Alikakos--Fusco \cite{AF2}.  Once \eqref{eq:maxp} is known, set $v=|u-a_i|^2$.  By \eqref{EL eq}, $v$ satisfies
\begin{equation*}
 \e^2\Delta v\geq cv,\quad\text{in }\Om.
\end{equation*}
\eqref{eq:exp decay} follows from standard comparison principle in elliptic theory. 
\end{proof}

\section{Minimizers on the unit disk with well-prepared boundary data}\label{sec:min unit disk}

Let $\Gamma=L_1\cup L_2$ be a cross with center $O_\G=(x_0,y_0)\in B_{\f1{10}}$. And suppose $L_1$, $L_2$ have the expressions \eqref{cross lines} with the half-angle of the $D_1$ sector satisfying \eqref{strict-cross-condition}. We focus on the restriction of $\G$ on the unit ball $B_1$. Denote by $P_1,P_3$ the intersections of $L_1$ with $\pa B_1$ such that $x_{P_1}>x_{P_3}$, and denote by $P_2,P_4$ the intersections of $L_2$ with $\pa B_1$ such that $x_{P_2}>x_{P_4}$. Set
\begin{equation*}
    l_i:=\dist(P_i,O_\G).
\end{equation*}
For convenience of calculation, throughout this section we translate the coordinates by $O_\G$. Thus the cross center is the origin and the disk is a shifted unit disk, with new center at $(-x_0,-y_0)$, denoted simply by $B$. Then we can write the coordinates of $\{P_i\}_{i=1}^4$ as 
\begin{equation}\label{cross-coordinates}
\begin{aligned}
 P_1&=l_1(\cos\phi,-\sin\phi),&
 P_2&=l_2(\cos\phi, \sin\phi),\\
 P_3&=l_3(-\cos\phi,\sin\phi),&
 P_4&=l_4(-\cos\phi,-\sin\phi).
\end{aligned}
\end{equation}

Set $B\setminus \Gamma =\cup_{i=1}^4 D_i$, where $D_i's$ are four sectors centered at the origin (now is also the center $O_\G$), labeled in counterclockwise order so that $D_1$ is the rightmost subdomain. Therefore, $\pa D_i\cap \pa B$ is the circular arc between $P_i$ and $P_{i+1}$, for $i=1,2,3,4$ ($P_5=P_1$).

We next prescribe the Dirichlet boundary data $g_{\e,\G}\in C^\infty(\pa B,\BR^2)$ associated with $\G$ for sufficiently small $\e$.  On the boundary arc from $P_i$ to $P_{i+1}$, after removing $\e$-neighborhoods of endpoints, set $g_{\e,\G}\equiv a_i$.  Denote the transition arc centered at $P_i$ by $G_i$, which has length $2\e$. Uniformly in $\e$, we assume 
\begin{equation}\label{est:trace}
 \|g_{\e,\G}\|_\infty+\e\|\partial_\tau g_{\e,\G}\|_\infty\leq C,
 \qquad
 \int_{\partial B}\left(\f{\e}{2}|\partial_\tau g_{\e,\G}|^2+\f{1}{\e}W(g_{\e,\G})\right)\,d\mathcal H^1\leq C.
\end{equation}

\begin{figure}[htbp!]
\centering
\begin{tikzpicture}[scale=2.45,line cap=round,line join=round]
  \coordinate (P1) at (-43:1);
  \coordinate (P2) at (47:1);
  \coordinate (P3) at (137:1);
  \coordinate (P4) at (227:1);
  \draw[black!45,line width=.5pt] (0,0) circle (1);

  \draw[red!50,line width=2.1pt,opacity=.5]
    (-38:1) arc[start angle=-38,end angle=42,radius=1];
  \draw[blue!50,line width=2.1pt,opacity=.5]
    (52:1) arc[start angle=52,end angle=132,radius=1];
  \draw[green!50,line width=2.1pt,opacity=.5]
    (142:1) arc[start angle=142,end angle=222,radius=1];
  \draw[purple!50,line width=2.1pt,opacity=.5]
    (232:1) arc[start angle=232,end angle=312,radius=1];

  \draw[purple,line width=4pt] (-48:1)
    arc[start angle=-48,end angle=-43,radius=1];
  \draw[red,line width=4pt] (-43:1)
    arc[start angle=-43,end angle=-38,radius=1];
  \draw[red,line width=4pt] (42:1)
    arc[start angle=42,end angle=47,radius=1];
  \draw[blue,line width=4pt] (47:1)
    arc[start angle=47,end angle=52,radius=1];
  \draw[blue,line width=4pt] (132:1)
    arc[start angle=132,end angle=137,radius=1];
  \draw[green,line width=4pt] (137:1)
    arc[start angle=137,end angle=142,radius=1];
  \draw[green,line width=4pt] (222:1)
    arc[start angle=222,end angle=227,radius=1];
  \draw[purple,line width=4pt] (227:1)
    arc[start angle=227,end angle=232,radius=1];

  \draw[black,thick] (P1)--(P3);
  \draw[black,thick] (P2)--(P4);
  \draw[dashed] (0,0)--(0.45,0);
  \draw (0.21,0) arc[start angle=0,end angle=48,radius=0.21];
    \node at ({0.32*cos(24)},{0.32*sin(24)}) {$\phi$};

  \foreach \p/\lab in {P1/P_1,P2/P_2,P3/P_3,P4/P_4}{
    \fill[black] (\p) circle (.024);
    \node[font=\scriptsize,fill=white,inner sep=.6pt]
      at ($ (0,0)!1.15!(\p) $) {\(\lab\)};
  }
  \node[text=red!80!black] at (.82,0) {\(a_1\)};
  \node[text=blue!80!black] at (0,.82) {\(a_2\)};
  \node[text=green!80!black] at (-.82,0) {\(a_3\)};
  \node[below] at (0,0) {\small$O_\G=(0,0)$};
  \node[font=\small,text=purple!80!black] at (0,-.82) {\(a_4\)};

  \node[anchor=west] (transitionlabel) at (1.08,-.5)
    {$G_i:$ transition arc of length $2\e$};
  \draw[->,thin] (transitionlabel.south west) -- (-40.5:1.03);
\end{tikzpicture}
\caption{The trace $g_{\e,\G}$ on $\pa B$.  It is constant on the four colored phase arcs and transits smoothly in thick $O(\e)$ transition arcs centered at $P_i$.}
\label{fig:trace}
\end{figure}
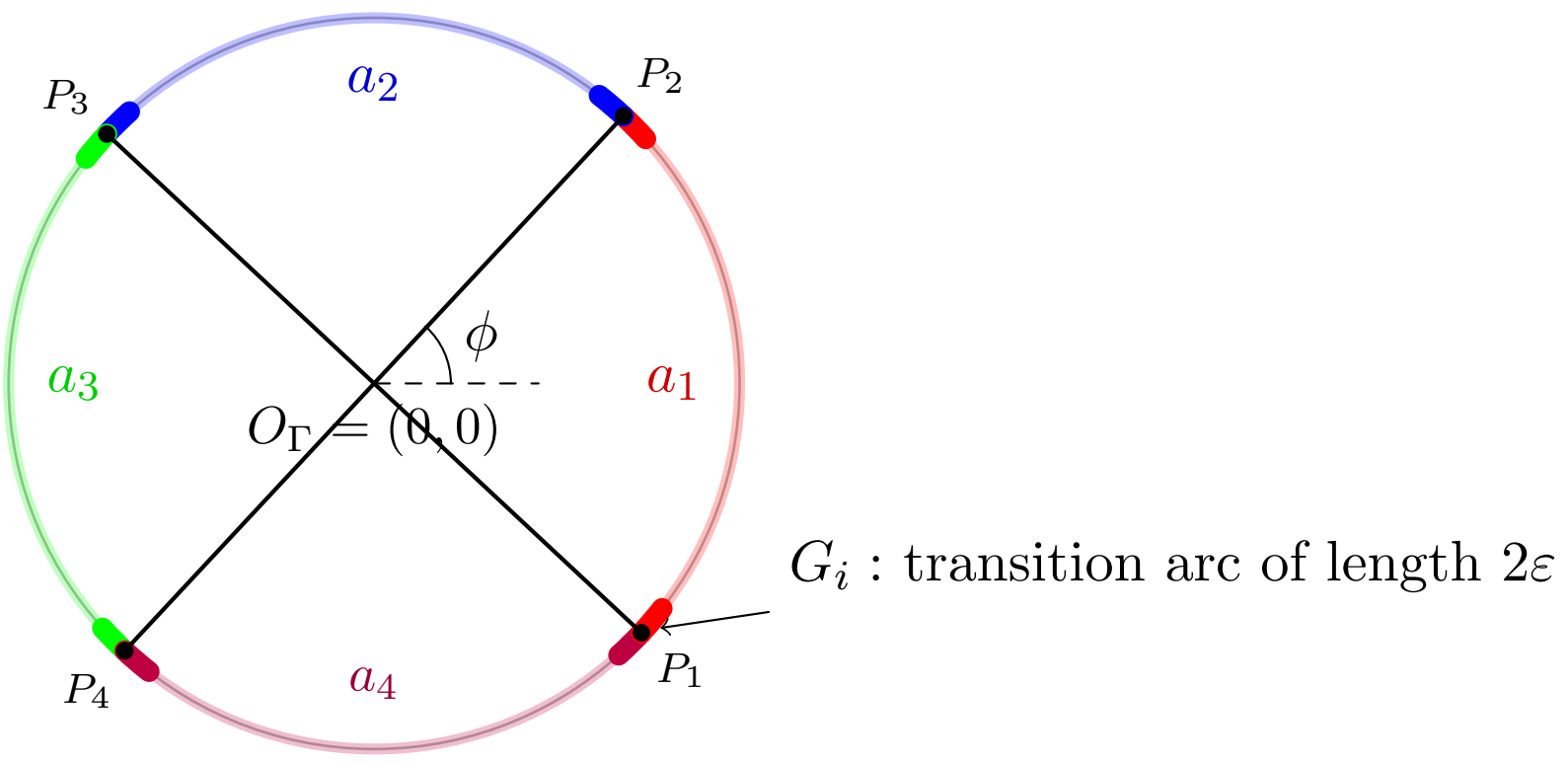

Let $u_\e$ minimize $J_\e(\cdot;B)$ subject to Dirichlet boundary condition $g_{\e,\G}$.  The existence of a minimizer follows from the direct method of calculus of variations. Using \eqref{radial coercivity} and standard theory of elliptic regularity, one can establish the $L^\infty$ bound for $u_\e$ and then bootstrap to $C^{2,\alpha}$ regularity. 

We start with the energy upper bound.

\begin{lemma}[Upper bound]\label{lem:upper bdd}
Let $\G$ and $B$ be as above. Suppose $u_\e$ is a minimizer with the trace $u_\e=g_{\e,\G}$ on $\pa B$ and $\e$ is sufficiently small, then
\begin{equation}\label{eq:upper bdd}
 J_\e(u_\e;B)\leq \sum_{i=1}^4l_i+C\e,
\end{equation}
where the constant $C=C(W)$ is independent of $\e$.
\end{lemma}

\begin{proof}
The upper bound follows from a standard energy competitor construction obtained by inserting one-dimensional transition profiles into a small neighborhood of prescribed minimal network; see for instance \cite[Appendix A]{alikakos2024triple} for detailed constructions.  We recall only the main idea.

Construct a competitor $w$ by setting $w\equiv a_i$ away from a thin neighborhood of the cross in the sector $D_i$.  In an small tube around each sharp interface between $D_i$ and $D_j$, replace the sharp jump between phases $a_i$ and $a_j$ by the corresponding truncated heteroclinic connection $U_{ij}$,  evaluated at signed distance from $\pa D_i\cap\pa D_j$ divided by $\e$. The four tubes meet in an $O(\e)$-neighborhood of the junction point $O_\G$, which we call the junction region. On the boundary of this $O(\e)$ region, the preceding construction gives a trace consisting of four transition arcs separated by four constant phase arcs. Then one can construct a Lipschitz extension, with Lipschitz constant bounded by $\frac{C}{\e}$, whose energy inside the junction region is $O(\e)$. Similarly, in an $O(\e)$ boundary layer of $\pa B$, smoothly interpolating between the interior construction and the prescribed trace $g_{\e,\G}$ yields an $O(\e)$ energy on the boundary collar.

Using equipartition and exponential decay of the heteroclinics in Lemma \ref{lem:heteroclinic}, the four tubes around $\G$ contribute $\sum_{i=1}^4l_i+O(\e)$.  Therefore
\begin{equation*}
 J_\e(w;B)\le\sum_{i=1}^4l_i+C\e
\end{equation*}
This completes the upper bound proof thanks to the minimality of $u_\e$.
\end{proof}

We define the following notations.

    \begin{align*}
        &A_i:=\{|u_\e(x)-a_i|\leq \e^{1/6}\},\quad i=1,2,3,4;\\
        &V_t:= \{(x,y): x=t\}\cap B, \quad \text{vertical slice in }B,\\
        &H_t:= \{(x,y): y=t\}\cap B, \quad \text{horizontal slice in }B;\\
        &j_v(t):= \int_{V_t} |\pa_y u_\e| \sqrt{2W(u_\e)}\,dy, \quad t\in(-1-x_0,1-x_0),\\
        &j_h(t):= \int_{H_t} |\pa_x u_\e| \sqrt{2W(u_\e)}\,dx, \quad t\in(-1-y_0,1-y_0);\\
        &\J_v(U):= \int_U |\pa_y u_\e|\sqrt{2W(u_\e)}\,dxdy,\quad \text{ for any measurable }U\subset B,\\
        &\J_h(U):= \int_U |\pa_x u_\e|\sqrt{2W(u_\e)}\,dxdy,\quad \text{ for any measurable }U\subset B.
    \end{align*}

Note that here the choices for the range of $t$ are to guarantee that $V_t$ or $H_t$ is nonempty, given that the center of $B$ is $(-x_0,-y_0)$ after the translation of the coordinates. Moreover, since $W(u_\e)\geq c\e^{1/3}$ on $B\setminus(\cup_{i=1}^4 A_i)$, Lemma \ref{lem:upper bdd} implies that 
\begin{equation}\label{est:measure of diffuse intf}
    |B\setminus(\cup_{i=1}^4 A_i) |\leq C\e^{2/3}. 
\end{equation}

We further define
\begin{align*}
N_v(t)&:=\mathbf{1}_{(0,l_1\cos\phi)}(t)+\mathbf{1}_{(0,l_2\cos\phi)}(t)+\mathbf{1}_{(-l_3\cos\phi,0)}(t)+\mathbf{1}_{(-l_4\cos\phi,0)}(t), \quad t\in(-1-x_0,1-x_0),\\
N_h(t)&:=\mathbf{1}_{(0,l_2\sin\phi)}(t)+\mathbf{1}_{(0,l_3\sin\phi)}(t)+\mathbf{1}_{(-l_1\sin\phi,0)}(t)+\mathbf{1}_{(-l_4\sin\phi,0)}(t),\quad t\in(-1-y_0,1-y_0).
\end{align*}
It's not hard to see that $N_v(t)$ (or $N_h(t)$) is the number of intersections of the vertical slice $V_t$ (or horizontal slice $H_t$) with $(\G\cap B)\setminus \{0\}$.

Next we try to define several types of ``bad slices", all of which will eventually be shown, either directly or with some effort, to have small measure. For a vertical slice $V_t$, we call it a ``Type-I" bad slice if it intersects with the transition arcs $G_i$ on $\pa B$, or it contains a subset of measure at least $\e^{1/3}$, on which $\min|u_\e-a_i|>\e^{1/6}$. More precisely, the set of ``Type-I" bad vertical slices is defined by
\begin{equation*}
    \mathcal{B}_1^v:=\{t\in (-1-x_0,1-x_0): \overline{V_t}\cap(\bigcup_{i=1}^4G_i)\neq \emptyset\ \text{ or }\  \ch(V_t\setminus (\bigcup_{i=1}^4 A_i))\geq \e^{1/3}\}.
\end{equation*}
Similarly, the set of ``Type-I" bad horizontal slices is defined by
\begin{equation*}
    \mathcal{B}_1^h:=\{t\in (-1-y_0,1-y_0): \overline{H_t}\cap(\bigcup_{i=1}^4G_i)\neq \emptyset\ \text{ or }\  \ch(H_t\setminus (\bigcup_{i=1}^4 A_i))\geq \e^{1/3}\}.
\end{equation*}
It follows immediately from the $O(\e)$ size of $G_i$ and \eqref{est:measure of diffuse intf} that 
\begin{equation*}
\ch(\mathcal{B}_1^h)+\ch(\mathcal{B}_1^v)\leq C\e^{1/3}.
\end{equation*}

Consider a vertical slice $V_t$ such that $N_v(t)=2$ and $t\notin \mathcal{B}_1^v$. By definition, $V_t$ intersects $\pa B$ at two points, where $g_{\e,\G}$ takes value $a_2$ and $a_4$ respectively, and passes by the region $D_3$ (if $t<0$) or $D_1$ (if $t>0$). We call it a ``Type-II" bad slice if $V_t$ avoids $A_3\cup A_1$. More precisely, set
\begin{equation*}
    \mathcal{B}_2^v:=\{t:N_v(t)=2,\ V_t\cap (A_1\cup A_3)=\emptyset\} \setminus \mathcal{B}_1^v. 
\end{equation*}
Similarly, horizontal ``Type-II" bad slice set is set to be
\begin{equation*}
    \mathcal{B}_2^h:=\{t:N_h(t)=2,\ H_t\cap (A_2\cup A_4)=\emptyset\} \setminus \mathcal{B}_1^h. 
\end{equation*}
Further define the essential suprema and infima of ``Type-II" bad sets:
\begin{align*}
   &\text{If }\ |\mathcal{B}_2^h|>0,\ \   y_i:= \essinf \B_2^h,\ \ y_s:= \esssup \B_2^h;\quad \text{otherwise }\ y_i=y_s=0.\\
   &\text{If }\ |\mathcal{B}_2^v|>0, \ \   x_i:= \essinf \B_2^v,\ \ x_s:= \esssup \B_2^v;\quad \text{otherwise }\ x_i=x_s=0.
\end{align*}

\begin{lemma}\label{lem:type-II}
Let $\B_2^h,\,\B_2^v,\,y_i,\,y_s,\,x_i,\,x_s$ be as above. There exists a constant $C=C(W)$ such that 
\begin{equation}\label{est:type-ii size}
    |\B_2^h|+y_s-y_i \leq C\e^{1/3},\qquad  |\B_2^v|+x_s-x_i \leq C\e^{1/3},
\end{equation}
and 
\begin{equation}\label{est:type-ii location}
|y_i|+|y_s|+|x_i|+|x_s|\leq C\e^{1/6}.    
\end{equation}
Moreover, the following estimates for horizontal and vertical energies hold:
\begin{align}
    \label{Jh lower bdd} & \J_h(B)\geq  \sin\phi (\sum_{i=1}^4l_i)-C\e^{1/3},\\
    \label{Jv lower bdd} &\J_v(B)\geq  \cos\phi (\sum_{i=1}^4l_i)-C\e^{1/3}.
\end{align}
\end{lemma}

\begin{proof}
    If $|\B_2^h|=|\B_2^v|=0$, then \eqref{est:type-ii size} and \eqref{est:type-ii location} hold trivially. We examine the horizontal slices $H_t$ for $t\in (-1-y_0,1-y_0)\setminus (\B_1^h\cup\B_2^h)$. When $N_h(t)=1$, $H_t$ intersects $\pa B$ at two points of adjacent phases (i.e. $a_3a_2$, $a_2a_1$, $a_3a_4$ or $a_4a_1$), which implies that at least one phase transition occurs along $H_t$. By \eqref{est:endpoint close to ai} we have
    \begin{equation*}
        j_h(t)\geq 1-C\e^{1/3}.
    \end{equation*}
    When $N_h(t)=2$, since $t\notin \B_2^h$, $H_t$ must intersect with $\pa B$ at an $a_3$ point and an $a_1$ point, and also pass by a point in $A_2$ or $A_4$. Thus, there are roughly at least two phase transitions on $H_t$, which yields 
    \begin{equation*}
        j_h(t)\geq 2-C\e^{1/3}.
    \end{equation*}
    As a result, we have 
    \begin{equation*}
        \J_h(B)=\int_{-1-y_0}^{1-y_0} j_h(t)\,dt\geq \int_{-1-y_0}^{1-y_0} N_h(t)\,dt-C\e^{1/3}=\sin\phi(\sum_{i=1}^4 l_i)-C\e^{1/3},
    \end{equation*}
    which gives \eqref{Jh lower bdd}. \eqref{Jv lower bdd} follows from the same calculation when $|\B_2^v|=0$.

Now we consider the situation when at least one of $|\B_2^v|,\,|\B_2^h|$ is nonzero.   Assume without loss of generality that $|\B_2^h|>0$, and consequently
    \begin{equation*}
        \rho:=y_s-y_i>0. 
    \end{equation*}
For convenience, we may assume $y_s,y_i\in \B_2^h$. Otherwise, we can take a sequence $y_{s,n}\to y_s$ and $y_{i,n}\to y_i$ such that $\{y_{i,n}\}\cup \{y_{s,n}\}\subset \B_2^h$. Then we can work with $y_{i,n},\,y_{s,n}$, derive uniform estimates and finally let $n\to\infty$. 

Now for every $t\in (y_i,y_s)\setminus \B_1^h$, $H_t$ joins $a_3$ with $a_1$, which implies 
\begin{equation}\label{est Jh yi<y<ys}
    \J_h( \{(x,y)\in B: y_i<y<y_s\})\geq s\rho-C\e^{1/3}.
\end{equation}
On the other hand, following the same argument as in the situation $|\B_2^h|=0$, for almost every $t\in \big((-1-y_0,y_i)\cup(y_s,1-y_0)\big)\setminus \B_1^h$, $j_h(t)$ will be bounded below by $N_h(t)-C\e^{1/3}$, which gives
\begin{align}
 \label{est:Jh, y>ys}   &\J_h(B\cap\{y>y_s\}) \geq \int_{y_s}^{1-y_0} N_h(t)\,dt-C\e^{1/3}=\sin\phi(l_2+l_3)-2y_s-C\e^{1/3},\\
 \label{est:Jh, y<yi}   &\J_h(B\cap\{y<y_i\}) \geq \int_{-1-y_0}^{y_i} N_h(t)\,dt-C\e^{1/3}=\sin\phi(l_1+l_4)+2y_i-C\e^{1/3}.
\end{align}

We now examine the vertical line $V_x$. On the upper subdomain $B\cap\{y>y_s\}$, for most of $x\in (-l_3\cos\phi,l_2\cos\phi)$ (minus a small set with measure less than $\e^{1/3}$), $V_x$ intersects with $\pa B$ at an $a_2$ point and with $H_{y_s}$ at an $a_3$ or $a_1$ point. We have 
\begin{equation}\label{est:Jv y>ys}
    \J_v(B\cap\{y>y_s\})\geq \cos\phi(l_2+l_3)-C\e^{1/3}.
\end{equation}
Similarly,
\begin{equation}\label{est:Jv y<yi}
    \J_v(B\cap\{y<y_i\})\geq \cos\phi (l_1+l_4)-C\e^{1/3}.
\end{equation}

For any $\kappa\in (0,1)$ and $U\subset B$, we have the following estimate
\begin{align*}
    J_\e(u_\e,U)&=\int_U\left\{\left(\frac{\e}{2}|\pa_x u_\e|^2+\frac{\kappa^2}{\e}W(u_\e)  \right)+\left(\frac{\e}{2}|\pa_y u_\e|^2+\frac{1-\kappa^2}{\e}W(u_\e)  \right)\right\}\,dxdy\\
    &\geq \kappa \J_h(U)+\sqrt{1-\kappa^2} \J_v(U). 
\end{align*}
By choosing $\kappa$ such that $\frac{\kappa}{\sqrt{1-\kappa^2}}=\f{\J_h(U)}{\J_v(U)}$, we obtain
\begin{equation}\label{ineq:Je Jv Jh}
   J_\e(u_\e,U)\geq \sqrt{|\J_v(U)|^2+|\J_h(U)|^2}. 
\end{equation}
Applying \eqref{ineq:Je Jv Jh} to $B\cap\{y>y_s\}$ and $B\cap \{y<y_i\}$, together with \eqref{est:Jh, y>ys}, \eqref{est:Jh, y<yi}, \eqref{est:Jv y>ys} and \eqref{est:Jv y<yi} implies
\begin{align}
    &J_\e(u_\e, B\cap\{y>y_s\}) \geq \sqrt{(\sin\phi(l_2+l_3)-2y_s)^2+(\cos^2\phi)(l_2+l_3)^2}-C\e^{1/3},\\
    &J_\e(u_\e, B\cap\{y<y_i\}) \geq \sqrt{(\sin\phi(l_1+l_4)+2y_i)^2+(\cos^2\phi)(l_1+l_4)^2}-C\e^{1/3}
\end{align}
Combining with \eqref{est Jh yi<y<ys} gives
\begin{equation}\label{est: Je lower 1}
\begin{split}
    J_\e(u_\e,B)&\geq \sqrt{(\sin\phi(l_2+l_3)-2y_s)^2+(\cos^2\phi)(l_2+l_3)^2}\\
    &\quad +\sqrt{(\sin\phi(l_1+l_4)+2y_i)^2+(\cos^2\phi)(l_1+l_4)^2}+s(y_s-y_i)-C\e^{1/3}.
\end{split}
\end{equation}
Define the function
\begin{equation*}
    f(t):= \sqrt{(\sin\phi(l_2+l_3)-2t)^2+(\cos^2\phi)(l_2+l_3)^2}+st.
\end{equation*}
Direct computation shows
\begin{align*}
    &f'(t)= \frac{4t-2\sin\phi(l_2+l_3)}{\sqrt{(\sin\phi(l_2+l_3)-2t)^2+(\cos^2\phi)(l_2+l_3)^2}}+s,\\
    &f'(0)=s-2\sin\phi>0,\\
    &f''(t)= \frac{4\cos^2\phi\,(l_2+l_3)^2}{\left( (\sin\phi(l_2+l_3)-2t)^2+(\cos^2\phi)(l_2+l_3)^2 \right)^{3/2}}>0.
\end{align*}
Thus by Taylor expansion we have 
\begin{equation*}
    f(t)\geq l_2+l_3+(s-2\sin\phi)t+ct^2,
\end{equation*}
for some positive constant $c>0$. 

Therefore, from \eqref{est: Je lower 1} we obtain
\begin{equation}\label{est: Je lower 2}
    J_\e(u_\e,B)\geq \sum_{i=1}^4 l_i+(s-2\sin\phi)(y_s-y_i)+ c(|y_s|^2+|y_i|^2)-C\e^{1/3}.
\end{equation}
Compare this with the energy upper bound \eqref{eq:upper bdd}, we have  
\begin{align*}
&\qquad (s-2\sin\phi)(y_s-y_i)+c(|y_s|^2+|y_i|^2)\leq C\e^{1/3},\\
&\Rightarrow\quad  
|y_s-y_i|\leq C\e^{1/3}, \quad |y_s|+|y_i|\leq C\e^{1/6}, 
\end{align*}
which gives \eqref{est:type-ii size} and \eqref{est:type-ii location} for $y_s,\,y_i$. The estimates for $x_s,\,x_i$ can be proved in the same manner when $|\B_2^v|>0$. Finally, the horizontal and vertical energy estimates \eqref{Jh lower bdd} and \eqref{Jv lower bdd} follow immediately from \eqref{est:Jh, y>ys}--\eqref{est:Jv y<yi}.

\end{proof}

\begin{corol}\label{corol: lower bound}
    Let $u_\e$ be a minimizer of $J_\e$ with the trace $u_\e=g_{\e,\G}$ on $\pa B$, there exists a constant $C=C(W)$ such that for sufficiently small $\e$, 
    \begin{equation}
        \label{eq:lower bdd}
        J_\e(u_\e,B)\geq \sum_{i=1}^4 l_i-C\e^{1/3}.
    \end{equation}
\end{corol}
\begin{proof}
    \eqref{eq:lower bdd} follows directly from \eqref{Jh lower bdd}, \eqref{Jv lower bdd} and \eqref{ineq:Je Jv Jh}.
\end{proof}

\begin{lemma}\label{lem: Jh Jv upper bdd}
    The following upper bounds for horizontal and vertical energies hold
   \begin{align}
    \label{Jh upp bdd} & \J_h(B)\leq  \sin\phi (\sum_{i=1}^4l_i)+C\e^{1/3},\\
    \label{Jv upp bdd} &\J_v(B)\leq  \cos\phi (\sum_{i=1}^4l_i)+C\e^{1/3}.
\end{align}
\end{lemma}
\begin{proof}
    By \eqref{ineq:Je Jv Jh}, \eqref{eq:upper bdd} and \eqref{Jv lower bdd}, we have 
    \begin{align*}
        &\J_h(B)\leq \sqrt{J_\e(u_\e,B)^2-\J_v(B)^2}\\
        &\qquad \leq \sqrt{(\sum_{i=1}^4l_i)^2+C\e-\cos^2\phi\,(\sum_{i=1}^4l_i)^2+C\e^{1/3} }\leq \sin\phi\sum_{i=1}^4l_i +C\e^{1/3},
    \end{align*}
    which implies \eqref{Jh upp bdd}. The same argument also applies to $\J_v(B)$.
\end{proof}

Our next key lemma shows the balanced relation between the horizontal energy $\J_h$ and vertical energy $\J_v$ on any measurable set $U\subset B$. 
\begin{lemma}\label{lem: balance h v}
    There is a constant $C=C(W)$ such that for every measurable $U\subset B$, 
    \begin{equation*}
        |\sin\phi\, \J_v(U)-\cos\phi\, \J_h(U)|\leq C\e^{1/6}.
    \end{equation*}
\end{lemma} 
\begin{proof}
    Fix an arbitrary measurable set $U\subset B$. Direct computation implies 
    \begin{align*}
         &|\sin\phi\, \J_v(U)-\cos\phi\, \J_h(U)|\\
         &=\left|  \int_{U}\sqrt{\f{2W(u_\e)}{\e}}\left(\sqrt{\e}|\pa_yu_\e|\sin\phi-\sqrt{\e}|\pa_xu_\e|\cos\phi\right)\,dz\right|\\
         &=\left|  \int_{U}\sqrt{\f{2W(u_\e)}{\e}}\left(\left(\sqrt{\e}|\pa_yu_\e|-\cos\phi \sqrt{\f{2W(u_\e)}{\e}}\right)\sin\phi-\left(\sqrt{\e}|\pa_xu_\e|-\sin\phi\sqrt{\f{2W(u_\e)}{\e}}\right)\cos\phi\right)\,dz\right|\\
         &\leq \left( \int_U \frac{2W(u_\e)}{\e} \right)^{1/2}\left( \int_U \left[\left( \sqrt{\e}|\pa_xu_\e|-\sin\phi\sqrt{\f{2W(u_\e)}{\e}} \right)^2+\left(\sqrt{\e}|\pa_yu_\e|-\cos\phi\sqrt{\f{2W(u_\e)}{\e}}\right)^2 \right]\,dz\right)^{1/2}\\
         &\leq \left( \int_B \frac{2W(u_\e)}{\e} \right)^{1/2}\left( \int_B \left[\left( \sqrt{\e}|\pa_xu_\e|-\sin\phi\sqrt{\f{2W(u_\e)}{\e}} \right)^2+\left(\sqrt{\e}|\pa_yu_\e|-\cos\phi\sqrt{\f{2W(u_\e)}{\e}}\right)^2 \right]\,dz\right)^{1/2}\\
         &\leq \sqrt{2J_\e(u_\e, B)} \sqrt{ 2(J_\e(u_\e,B)-\sin\phi\,\J_h(B)-\cos\phi\,\J_v(B)) }\\
         &\leq C\e^{1/6}.
    \end{align*}
    Here we used \eqref{eq:upper bdd}, \eqref{Jh lower bdd} and \eqref{Jv lower bdd} in the last step.
\end{proof}

Next we define the last type of bad slices, which essentially contains all the slices that cost more energy than the minimal energy required by the visited phases on it. More precisely, we set the ``Type-III" bad horizontal  set and bad vertical set as
\begin{align*}
    &\B_3^h:= \{y\in (-1-y_0,1-y_0)\setminus (\B_1^h\cup \B_2^h): j_h(y)\geq N_h(y)+\e^{1/6} \},\\
    &\B_3^v:= \{x\in (-1-x_0,1-x_0)\setminus (\B_1^v\cup \B_2^v): j_v(x)\geq N_v(x)+\e^{1/6} \}.
\end{align*}

\begin{lemma}\label{lem:type-iii}
    There is a constant $C=C(W)$ such that 
    \begin{equation*}
        |\B_3^h|+|\B_3^v|\leq C\e^{1/6}.
    \end{equation*}
\end{lemma}

\begin{proof}
It suffices to consider only $\B_3^h$. The argument naturally extends to $\B_3^v$ without any change. For any $y\in (-1-y_0,1-y_0)\setminus (\B_1^h\cup \B_2^h)$, it holds that
\begin{align*}
    &N_h(y)=1 \ \Rightarrow\ j_h(y)\geq 1-C\e^{1/3},\\
    &N_h(y)=2 \ \Rightarrow\ j_h(y)\geq 2-C\e^{1/3}.
\end{align*}
We have 
\begin{align*}
    \J_h(B)&\geq \int_{(-1-y_0,1-y_0)\setminus (\cup_{i=1}^3 \B_i^h) }j_h(y)\,dy+ \int_{\B_3^h} j_h(y)\,dy\\
    \geq & \sin\phi(\sum_{i=1}^4l_i)+\e^{1/6}|\B_3^h|-C\e^{1/3},
\end{align*}
This, together with \eqref{Jh upp bdd}, implies
\begin{equation*}
    \e^{1/6} |\B_3^h|\leq C\e^{1/3}\ \Rightarrow \ |\B_3^h|\leq C\e^{1/6}.
\end{equation*}
\end{proof}

\begin{lemma}\label{lem: Jh, Jv subset} For any measurable set $Y\subset (-1-y_0,1-y_0)$,
\begin{equation}\label{est:Jh subset}
    \int_{Y} j_h(y)\,dy\leq \int_{Y} N_h(y)\,dy+C\e^{1/3}.
\end{equation}
Similarly, for any measurable set $X\subset (-1-x_0,1-x_0)$,
\begin{equation}\label{est:Jv subset}
    \int_{X} j_v(x)\,dx\leq \int_{X} N_v(x)\,dx+C\e^{1/3}.
\end{equation}
    
\end{lemma}

\begin{proof}
We have
\begin{align*}
    &\int_Y j_h(y)\,dy \leq \J_h(B)- \int_{(-1-y_0,1-y_0)\setminus (Y \cup \B_1^h\cup \B_2^h)  } j_h(y)\,dy\\
    &\qquad \leq \sin\phi(\sum_{i=1}^4l_i)+C\e^{1/3}- \int _{(-1-y_0,1-y_0)\setminus (Y \cup \B_1^h\cup \B_2^h)  } \left(N_h(y)-C\e^{1/3}\right)\,dy\\
    &\qquad=\int_{-1-y_0}^{1-y_0} N_h(y)\,dy+C\e^{1/3}- \int _{(-1-y_0,1-y_0)\setminus (Y \cup \B_1^h\cup \B_2^h)  } \left(N_h(y)-C\e^{1/3}\right)\,dy\\
    &\qquad \leq \int_{Y} N_h(y)\,dy +2|\B_1^h\cup\B_2^h|+C\e^{1/3}\\
    &\qquad \leq \int_YN_h(y)\,dy+C\e^{1/3},
\end{align*}
which proves \eqref{est:Jh subset}. \eqref{est:Jv subset} follows by the same calculation. 
\end{proof}

Now we are ready to present the main result of this section, which indicates the closeness of the diffuse interface to the sharp interface $\G$. 

\begin{proposition}[Localization of diffuse interface]\label{prop: loc diff int}
Let $\G,\,B,\,\{D_i\}_{i=1}^4,\,g_{\e,\G}$ be defined as in the beginning of Section \ref{sec:min unit disk}. Let $u_\e$ be a minimizer of $J_\e(\cdot,B)$ with the trace $u_\e=g_{\e,\G}$ on $\pa B$. For any $\gamma\in (0, \d_W)$, there exists a constant $C(\g,W)$ such that 
\begin{equation}\label{est: loc diff int}
    I_\e^\g\subset N_{C\e^{1/12}}(\Gamma):=\{x\in B: \mathrm{dist}(x, \G)\leq C\e^{1/12}\},
\end{equation}
where 
$$
I_\e^\g:=\{x\in B: \min_{1\leq i\leq 4}|u_\e(x)-a_i|>\g\}.
$$
Moreover, for every $x\in D_i$,
\begin{equation}\label{est: exponential}
    |u_\e(x)-a_i|\leq C\exp\left(-c\f{\max\{\dist(x,\pa D_i)-C\e^{1/12}, 0\}}{\e}\right),
\end{equation}
where $c,C$ are positive constants independent of $\e$.
\end{proposition}

The proof of Proposition \ref{prop: loc diff int} is based on the following key lemma.

\begin{lemma}\label{lem:4 phase radii}
    There are four pairwise disjoint line segments $\Sigma_1,\,\Sigma_2,\,\Sigma_3,\,\Sigma_4$ and constants $c,\,C$ depending only on $W$, such that the following properties hold: 
    \begin{enumerate}
        \item On $\Sigma_i$, $|u_\e-a_i|\leq c\e^{1/12}$.
        \item There exist $y_1,x_2,y_3,x_4$ such that 
        \begin{equation}\label{property of Sigma i}
        \begin{split}
        &\Sigma_1\subset D_1\cap H_{y_1},\qquad\Sigma_3\subset D_3\cap H_{y_3},\\
        &\Sigma_2\subset D_2\cap V_{x_2},\qquad \Sigma_4\subset D_4\cap V_{x_4} ;\\
        & |y_1|+|y_3|+|x_2|+|x_4| \leq C\e^{1/6}.
        \end{split}
        \end{equation}
        \item One end of $\Sigma_i$ belongs to the boundary arc $\arc{P_iP_{i+1}}$, the other end belongs to $B_{C\e^{1/6}}$. 
    \end{enumerate}
\end{lemma}

\begin{proof}[Proof of Lemma \ref{lem:4 phase radii}]  We construct $\Sigma_2$. Other segments can be constructed in the same way. 

Choose a constant $K$ such that 
\begin{equation}\label{choice K}
    |\B_1^h\cup\B_2^h\cup \B_3^h|+|\B_1^v\cup\B_2^v\cup\B_3^v|\leq \frac{K}{10(1+\tan\phi)}\e^{1/6}.
\end{equation}
In the subsequent argument, we may enlarge $K$ as needed. Then there exist $y_*\in (K\e^{1/6},\f{11}{10}K\e^{1/6})$ and $s<t$ satisfying the following properties:
\begin{enumerate}
\itemsep5pt
    \item[(i)] $y_* \notin \B_i^h$ for $i=1,2,3.$
    \item[(ii)] $A_3\cap H_{y_*} \subset \{(x,y_*):x<s\}$, $A_1\cap H_{y_*}\subset\{(x,y_*):x>t\}$.
    \item[(iii)] There is $i=2$ or $4$ such that $A_i\cap H_{y_*}\neq\emptyset$ and $A_i\cap H_{y_*}\subset\{(x,y_*):x\in(s,t)\}$.
    \item[(iv)] For every $i$ such that $A_i\cap H_{y_*}\neq \emptyset$, and for every point $z$ lying between the leftmost and the rightmost points of $\overline{A_i\cap H_{y_*}}$, it holds that $|u_\e(z)-a_i|\leq c\e^{1/12}$, for some constant $c=c(W)$. 
\end{enumerate}
(i) follows from our choice of constant $K$. Since $y_*$ is $O(\e^{1/6})$-close to $0$ and $y_*\notin \B_2^h$, $H_{y_*}$ must contain points in $A_3$, $A_1$, and either $A_2$ or $A_4$. We next justify property (iv). Let $(x_l,y_*)$, $(x_r,y_*)$ be the left and right most points of $\overline{A_i\cap H_{y_*}}$ respectively. Then for any $x\in [x_l,x_r]$, $(x,y_*)$ cannot belong to $A_j$ for some $j\neq i$. This is because $y_*\notin \B_3^h$, which means there are at most two transitions from $A_i$ to $A_j$ for $i\neq j$ if we traverse $H_{y_*}$ from left to right. Therefore, the phase transitions must occur outside the interval $[x_l,x_r]$ and by definition of $\B_3^h$ we have 
\begin{equation*}
    \int_{x_l}^{x_r} |\pa_x u_\e|\sqrt{2W(u_\e)}\,dx\leq C\e^{1/6}.
\end{equation*}
Since $(x_l,y_*),\,(x_r,y_*)\in A_i$, we conclude from \eqref{dW near well} that there exists $c=c(W)$ such that $|u_\e(x,y_*)-a_i|\leq c\e^{1/12}$ for every $x\in[x_l,x_r]$, which verifies property (iv). From (iv), properties (ii) and (iii) follow immediately since the three subsets $A_i\cap H_{y_*}$ are pairwise strictly separated in order on $H_{y_*}$. 

Next we claim that the middle phase on $H_{y_*}$ is $a_2$. Otherwise if the middle phase is $a_4$, then since $y_*\notin \B_1^h$, we have that except for a set of measure $\e^{1/3}$, $(x,y_*)\in A_3\cup A_4\cup A_1$, which implies
\begin{equation*}
    \J_v(\{y>y_*\}\cap B) \geq \cos\phi\,(l_2+l_3)-C\e^{1/3},
\end{equation*}
where $C$ only depends on $W$. By Lemma \ref{lem: balance h v}, we have 
\begin{equation}\label{est1: mid phase a2}
    \J_h(\{y>y_*\}\cap B)\geq \sin\phi(l_2+l_3)-C\e^{1/6}.
\end{equation}
On the other hand, \eqref{est:Jh subset} implies
\begin{equation}\label{est2:mid phase a2}
    \J_h(\{y>y_*\}\cap B)\leq \int_{y_*}^{1-y_0} N_h(y)\,dy+C\e^{1/3}\leq \sin\phi\,(l_2+l_3)-2y_*+C\e^{1/3},
\end{equation}
which yields a contradiction with \eqref{est1: mid phase a2} upon choosing a larger $K$ when necessary. 

Set
\begin{equation*}
    l:=\min\{x: (x,y_*)\in A_2\},\quad r:=\max\{x:(x,y_*)\in A_2\}.
\end{equation*}
Property (iv) guarantees that for every $x\in(l,r)$, $|u_\e(x,y_*)-a_2|\leq c\e^{1/12}$. Then the vertical phase-mismatch on $H_{y_*}$ and $\pa B\cap\{y>y_*\}$ implies
\begin{equation*}
    \J_v(\{y>y_*\}\cap B)\geq \cos\phi\,(l_2+l_3)-(r-l)-C\e^{1/3}.
\end{equation*}
Combining this with Lemma \ref{lem: balance h v} and \eqref{est2:mid phase a2} gives 
\begin{align*}
    &\cos\phi(l_2+l_3)-(r-l)-C\e^{1/3}\leq \cos\phi(l_2+l_3)-\f{2y_*}{\tan\phi}+C\e^{1/6}
    \\
    &\Rightarrow r-l\geq \f{2y_*}{\tan\phi}-C\e^{1/6}>\frac{19K}{10\tan\phi}\e^{1/6},
\end{align*}
after choosing larger $K$ if necessary. 

Now we locate the interval $(l,r)$. By \eqref{choice K}, there exists $x_*\in(l,l+\frac{K}{10\tan\phi}\e^{1/6})$ such that $(x_*,y_*)\in A_2$ and $x_*\notin \B^v_i$ for any $i\in\{1,2,3\}$. Set 
\begin{equation*}
\Om^+:=B\cap \{y>y_*,x>x_*\}, \qquad \Om^-:=B\cap \{y>y_*,x<x_*\}    
\end{equation*}
On the vertical segment $V_{x_*}\cap\{y>y_*\}$, both endpoints belong to $A_2$, then we may argue as in the derivation of property (iv) when selecting $y_*$ to conclude that 
\begin{equation*}
    |u_\e-a_2|\leq c\e^{1/12}, \quad \text{ on }\ V_{x_*}\cap\{y>y_*\}.
\end{equation*}
Therefore, 
\begin{equation*}
    \J_h(\Om^+)\geq (\sin\phi\, l_2-y_*)(1-C\e^{1/6}),
\end{equation*}
which further implies 
\begin{equation}
\label{loc1 r l}
\begin{split}
    &\J_h(\Om^-) =\J_h(B\cap\{y>y_*\})-\J_h(\Om^+)\\
    &\quad \leq \sin\phi(l_2+l_3)-2y_*+C\e^{1/3}-\J_h(\Om^+)\leq \sin\phi\,l_3-y_*+C\e^{1/6},
    \end{split}
\end{equation}
thanks to \eqref{est:Jh subset}. On the other hand, apart from a set of measure up to $\frac{K}{10\tan\phi}\e^{1/6}+\e^{1/3}$, $V_x$ will connect an $a_2$ boundary point to $(x,y_*)\in A_3$ when $x\in(-l_3\cos\phi, x_*)$. This implies
\begin{equation}\label{loc2 r l}
    \J_v(\Om^-)\geq (x_*+l_3\cos\phi-\frac{K}{10\tan\phi}\e^{1/6}-\e^{1/3})(1-C\e^{1/3})\geq \cos\phi\,l_3+l-\frac{K}{5\tan\phi}\e^{1/6}.
\end{equation}
Combining \eqref{loc1 r l}, \eqref{loc2 r l} and Lemma \ref{lem: balance h v} yields
\begin{equation*}
    (l-\frac{K}{5\tan\phi}\e^{1/6})\sin\phi\leq (-y_*+C\e^{1/6})\cos\phi+C\e^{1/6}\quad \Rightarrow\quad l\leq -\frac{3K}{5\tan\phi}\e^{1/6}. 
\end{equation*}
Similarly, we can establish
\begin{equation*}
    r\geq \frac{3K}{5\tan\phi}\e^{1/6}. 
\end{equation*}

Consequently, we can choose $x_2\in(l,r)$ such that 
\begin{equation*}
    (x_2,y_*)\in A_2,\quad |x_2|\leq \frac{K}{2\tan\phi}\e^{1/6},\quad |u-a_2|\leq c\e^{1/12} \ \text{  on }\ \Sigma_2:=V_{x_2}\cap\{y>y_*\}. 
\end{equation*}
It is straightforward to verify that $x_2$ and $\Sigma_2$ satisfy all properties required by Lemma \ref{lem:4 phase radii}. The proof is complete. 
\end{proof}

\begin{proof}[Proof of Proposition \ref{prop: loc diff int}]
Fix $c_0,\,C_0$ as the constants in the statement of Lemma \ref{lem:4 phase radii}. For each $i\in\{1,2,3,4\}$, we pick a point $Q_i\in \Sigma_i$, such that the distance of $Q_i$ to the interior endpoint of $\Sigma_i$ is $4C_0\e^{1/6}$. Let $R_i$ denote $\Sigma_i\cap \pa B$. Then we denote by $\Om_{i\,i+1}$ the region enclosed by the segments $Q_iR_i$, $Q_iQ_{i+1}$, $Q_{i+1}R_{i+1}$ and the arc $\arc{R_iR_{i+1}}$. From the definition, it is easy to see that $\Om_{12},\,\Om_{23},\,\Om_{34},\,\Om_{41}$ and the quadrilateral $Q_1Q_2Q_3Q_4$ form a complete partition of $B$.

Set
\begin{align*}
    &\tau_{12}:=(\cos\phi,\sin\phi),\quad\tau_{23}:=(-\cos\phi,\sin\phi),\quad \tau_{34}:=(-\cos\phi,-\sin\phi),\quad \tau_{41}:=(\cos\phi,-\sin\phi);\\
    &\nu_{12}:=(-\sin\phi,\cos\phi),\quad\nu_{23}:=(-\sin\phi,-\cos\phi),\quad \nu_{34}:=(\sin\phi,-\cos\phi),\quad \nu_{41}:=(\sin\phi,\cos\phi).
\end{align*}
Then $\tau_{i,i+1}$ represents the unit direction of $O_\G P_{i+1}$ and $\nu_{i,i+1}$ is the unit normal of $O_\G P_{i+1}$ pointing from $D_i$ to $D_{i+1}$. 

For any unit vector $\xi$ and measurable $U\subset B$, define the energy in the direction $\xi$ as
\begin{equation*}
    \J_{\xi}(U):=\int_U |\pa_\xi u_\e|\sqrt{2W(u_\e)} \,dz
\end{equation*}
Then the inequality \eqref{ineq:Je Jv Jh} can be generalized to 
\begin{equation}
    \label{ineq:Jxi Jxip}
    J_\e(u_\e,U)\geq \sqrt{|\J_\xi(U)|^2+|\J_{\xi^\perp}(U)|^2}.
\end{equation}

We first show that there exists $C_1(W)$ such that on every $\Om_{i\,i+1}$,
\begin{equation}\label{est:lower bdd Om_ii+1}
    \J_{\nu_{i\,i+1}}(\Om_{i\,i+1})\geq l_{i+1}-C_1\e^{1/6}.
\end{equation}
We work with $i=1$. For $t\in(5C_0\e^{1/6},l_2-\e)$, define 
\begin{equation*}
    L^\nu_t:= (t\cos\phi,t\sin\phi)+\nu_{12} \mathbb{R},
\end{equation*}
as the straight line orthogonal to $O_\G P_2$, whose distance to the origin is $t$. When $t\in (5C_0\e^{1/6},l_2-\e)$, $L^\nu_t$ intersects with $\pa \Om_{12}$ at two points $E_t$ and $F_t$, where $E_t$ lies on $Q_1R_1\cup \arc{R_1P_2}$ and $F_t$ lies on $Q_2R_2\cup \arc{P_2R_2}$. To see this, it suffices to notice that $L_t^\nu\cap Q_1Q_2=\emptyset$, since $|Q_i-O_\G|<5C_0\e^{1/6}$. We have 
\begin{equation*}
    |u_\e(E_t)-a_1|\leq c\e^{1/12},\qquad |u_\e(F_t)-a_2|\leq c\e^{1/12}.
\end{equation*}
Therefore,
\begin{equation*}
    \J_{\nu_{12}}(\Om_{12})\geq (l_2-5C_0\e^{1/6}-\e)(1-C\e^{1/6})\geq l_2-C_1\e^{1/6}.
\end{equation*}
We also obtain the upper bound for $J_\e(u_\e,\Om_{i\,i+1})$ :
\begin{equation}\label{upp bdd: J quad}
\begin{split}
    &J_\e(u_\e,\Om_{i\,i+1})\leq J_\e(u_\e,B)-\sum_{j\neq i} \J_{\nu_{j\,j+1}}(\Om_{j\,j+1})\\
    &\qquad \leq \sum_{i=1}^4 l_i+C\e- \left(\sum_{j\neq i}(l_{j+1} -C_1\e^{1/6}) \right)\leq l_{i+1}+3C_1\e^{1/6}.
\end{split}
\end{equation}

Combining \eqref{ineq:Jxi Jxip}, \eqref{est:lower bdd Om_ii+1} and \eqref{upp bdd: J quad}, we have 
\begin{equation}\label{est: quad par}
        \J_{\tau_{i\,i+1}} (\Om_{i\,i+1})\leq \sqrt{ (l_{i+1}+3C_1\e^{1/6})^2- (l_{i+1}-C_1\e^{1/6})^2}\leq C_2\e^{1/12},
\end{equation}
for $C_2$ only depending on $W$. 

For every point $z\in Q_2R_2\subset \Sigma_2$, let $\mathfrak{L}^{12}_z$ and $\mathfrak{L}^{23}_z$ denote the two segments joining $z$ to $\arc{P_2P_3}$ and parallel to direction $\tau_{12}$, $\tau_{23}$ respectively, i.e.
\begin{equation*}
    \mathfrak{L}^{12}_z:= \{z+\tau_{12}\mathbb{R}^+\}\cap B,\qquad \mathfrak{L}^{23}_z:= \{z+\tau_{23}\mathbb{R}^+\}\cap B.
\end{equation*}
We claim that there exists a constant $C_3(\g,W)$ and $z\in Q_2R_2$ such that 
\begin{equation}\label{choose z}
    \dist(z,Q_2)\leq C_3\e^{1/12},\quad |u_\e-a_2|\leq \g,\ \text{ on }\  \mathfrak{L}^{12}_z\cup\mathfrak{L}^{23}_z.
\end{equation}
\emph{Proof of the claim: }If there exists $z_1\in \mathfrak{L}^{12}_z$ satisfying $|u_\e(z_1)-a_2|>\g$, then 
\begin{equation*}
    \int_{\mathfrak{L}^{12}_z} |\pa_{\tau_{12}} u_\e|\sqrt{2W(u_\e)}\,d\ch\geq C\g^2.
\end{equation*}
Then \eqref{est: quad par} gives the $\ch$ measure of such $z$ is bounded by $\frac{C_2\e^{1/12}}{C\g^2}=C(\g,W)\e^{1/12}$. Similarly, the $\ch$ measure for the set of $z$, where $\mathfrak{L}^{23}_z$ contains a point of distance greater than $\g$ to $a_2$, is also controlled by $C(\g,W)\e^{1/12}$. Therefore, in the $3C(\g,W)\e^{1/12}$ neighborhood of $Q_2$ on $Q_2R_2$, one can find a point $z$ such that $|u_\e-a_2|\leq \g$ on $\mathfrak{L}^{12}_z\cup \mathfrak{L}^{23}_z$. The claim is proved by setting $C_3=3C(\g,W)$. 

Fix $z\in Q_2R_2$ satisfying \eqref{choose z}. Denote by $\mathcal{D}_{2}$ the region enclosed by $\mathfrak{L}^{12}_z$, $\mathfrak{L}^{23}_z$ and $\arc{P_2P_3}$. $\mathcal{D}_2$ is obtained from $D_2$ by removing a collar of width $C(\g,W)\e^{1/12}$ along $\pa D_2\cap \G$. We have 
\begin{equation*}
    |u_\e-a_2|\leq \g,\quad \text{ on }\pa \mathcal{D}_2.
\end{equation*}
Thus, we can invoke Proposition \ref{prop:maxp} to obtain
\begin{equation*}
    I_\e^\g\cap \mathcal{D}_2=\emptyset,
\end{equation*}
and 
\begin{equation*}
    |u_\e-a_2|\leq C\exp\left(-\f{c\dist(\cdot,\pa \mathcal{D}_2)}{\e}\right), \quad \text{in }\mathcal{D}_2. 
\end{equation*}

Similarly, one can construct $\mathcal{D}_i$ by removing a $C(\g,W)\e^{1/12}$-width-collar near $\G$ from $D_i$ for $i=1,3,4$, such that $u_\e$ is uniformly $\g$-close to $a_i$ in $\mathcal{D}_i$. Then
\begin{equation*}
    I_\e^\g\subset B\setminus (\bigcup_{i=1}^4 \mathcal{D}_i) \subset N_{C(\g,W)\e^{1/12}}(\G),
\end{equation*}
which leads to \eqref{est: loc diff int}. \eqref{est: exponential} follows directly from \eqref{eq:exp decay} applying to all $\mathcal{D}_i$. The proof of Proposition \ref{prop: loc diff int} is complete.
\end{proof}

\section{Perturbed trace and good radius}\label{sec:good radius}

The energy estimate and localization of diffuse interface in Section \ref{sec:min unit disk} are proved under the boundary data $g_{\e,\G}$ satisfying \eqref{est:trace}. Indeed, we can expand the class of admissible boundary data which still preserves the preceding fine structures of $u_\e$. We begin by describing the conditions that the expanded class of boundary data should satisfy.

Fix constants $\bar{\delta}>0$, $K\gg 1$, such that 
\begin{equation*}
    K\bar{\delta}< \frac{\delta_W}{4},\qquad c_W (K\bar{\delta})^2-C_W\bar{\delta}^2\geq 8C\bar{\delta}^2,
\end{equation*}
where $c_W,\,C_W$ are constants in \eqref{dW near well} and $C$ is the constant in \eqref{est:endpoint close to ai}. Also we need to choose another constant $C_{g}$ defined by
\begin{equation*}
    C_g:=\frac{5}{c\bar\delta^2},
\end{equation*}
where $c$ is the constant in \eqref{near well} such that $W(p)\geq c|p-a_i|^2$ for $p\in B_{2\delta_W}(a_i)$. Therefore $C_g$ depends on $W$ and the already chosen $\bar\delta$, which further implies it only depends on $W$. The reason for choosing such $K$ and $C_g$ will be made clear in the proof of Proposition \ref{prop: good radius}.

A function $g\in C^2 (\pa B_1;\BR^2)$ is called a \emph{good} boundary trace if the following hold

\begin{equation}\label{good partition}
\pa B_1=\left(\bigcup_{i=1}^4 G_i\right)\bigcup\left(\bigcup_{i=1}^4 \g_i\right),\qquad \text{all }G_i,\,\g_i\text{ are pairwise disjoint};
\end{equation}
\begin{equation}\label{good gap}
    \ch(G_i)\leq C_g\e,\qquad \min_{j=1,2,3,4}\ |g-a_j|\geq \bar\delta \ \ \text{ on }G_i,\qquad \forall\,i;
\end{equation}
\begin{equation}\label{good phase arc}
 |g-a_i|\leq K\bar\delta\ \ \text{ on }\g_i,\qquad \forall\,i;
\end{equation}
\begin{equation}\label{good energy}
        \int_{\pa B_1} \left(\frac{\e}{2}|\pa_\tau g|^2+\frac{1}{\e}W(g)\right)\,d\ch\leq 4+C\bar\delta^2,  \quad \text{ for constant }C \text{ in (\eqref{est:endpoint close to ai})}.
\end{equation}
Here the four small phase transition arcs $G_1,\,G_2,\,G_3,\,G_4$ occur counterclockwise. For every $i$, $\g_i$ is the $a_i$ phase arc between $G_i$ and $G_{i+1}$. We again let $P_i$ denote the midpoint of the transition gap $G_i$, define the associated cross and its center by
\begin{equation*}\label{good cross-definition}
 \Gamma=P_1P_3\cup P_2P_4\qquad O_\Gamma=P_1P_3\cap P_2P_4.
\end{equation*}
The following assumption on the cross angle ensures the minimality of $\G$:
\begin{equation}\label{good cross angle}
    \angle P_iO_\G P_{i+1}\in \bigg(2\arccos(\f{s}{2}),\ 2\arcsin(\f{s}{2})\bigg),\ \quad \text{for every }i.
\end{equation}
We finally make the following extra assumption on the cross center $O_\G$:
\begin{equation}\label{good center}
    |O_\G|\leq \f{1}{10}.
\end{equation}

\begin{lemma}\label{lem:fine prop good tr}
Let $g$ be a good trace satisfying \eqref{good partition}--\eqref{good center}, and let $u_\e$ minimize $J_\e(\cdot,B_1)$ subject to Dirichlet boundary condition $u_\e=g$ on $\pa B_1$. Then for $\e$ sufficiently small, we have the following energy bound:
\begin{equation}\label{ene bdd 3}
  |P_1P_3|+|P_2P_4|-C\e^{\f13}\leq J_\e(u_\e,B_1)\leq |P_1P_3|+|P_2P_4|+C\e.
  \end{equation}
  
Moreover, for every small $\g<\f{\delta_{W}}{2}$, we have the localization of the diffuse interface inside the annulus:
\begin{equation}\label{annulus loc dif int}
    \left(I_\e^\g\cap A_{\f14,\f34}(O_\G)\right)\  \subset\  N_{C(\g,W)\e^{1/12}}(\G\cap A_{\f14,\f34}(O_\G)),
\end{equation}
where $A_{r,R}(z)$ denotes the annulus $B_{R}(z)\setminus \overline{B_r}(z)$. 
\end{lemma}

\begin{proof}
For the cross $\G$ associated with the good boundary data $g$, we set for every $i\in\{1,2,3,4\}$
\begin{align*}
    \tilde{G}_i&:= \text{small arc centered at }P_i\text{ with length }2\e,\quad \tilde{G}_i\subset G_i,\\
    \tilde{\g}_i&:= \text{phase arc between }\tilde{G}_i,\,\tilde{G}_{i+1}, \quad \g_i\subset \tilde{\g}_i.
\end{align*}
Define the corresponding exact boundary data 
\begin{equation*}
    h=\begin{cases}
        a_i, &\text{ on }\tilde{\g}_i,\\
        \text{smoothly connect }a_{i-1},\,a_i,& \text{ on }\tilde{G}_i.
    \end{cases}
\end{equation*}
For the upper bound, construct the following energy competitor on $B_1$ in polar coordinates:
\begin{equation*}
    v_\e(re^{i\theta}):=\begin{cases}
        g(\theta), &\text{ on }\pa B_1,\\
        \f{1-r}{\e} h(\theta)+\f{r-1+\e}{\e} g(\theta),& r\in(1-\e,1),\\
        \text{minimizer of }J_\e(\cdot,B_{1-\e}) \text{ with trace }h, & \text{ on }B_{1-\e}.
    \end{cases}
\end{equation*}

Direct computation gives
\begin{align*}
 J_\e(v_\e;A_{1-\e,1})&\le C\int_{\mathbb S^1}|h-g|^2\,d\theta                        \\
 &\quad+C\e^2\int_{\mathbb S^1}\bigl(|\pa_\theta h|^2+|\pa_\theta g|^2\bigr)\,d\theta          \\
 &\quad+C\int_0^1\!\int_{\mathbb S^1} W\bigl((1-t)h+tg\bigr)\,d\theta \,dt.        \end{align*}
For the first term, we have 
\begin{align*}
 &\int_{\mathbb S^1}|h-g|^2\,d\theta  \leq \sum_i\int_{G_i}|h-g|^2\,d\theta+\sum_i \int_{\g_i}|g-a_i|^2\,d\theta\\
 &\qquad\quad \leq C\sum_i|G_i| +\int_{\mathbb S^1} W(g)\,d\theta \leq 4C\cdot C_g\e+C\e \leq C(W)\e.
\end{align*}
Similarly, the third term $\int_0^1\!\int_{\mathbb S^1} W\bigl((1-t)h+tg\bigr)\,d\theta \,dt$ can also be bounded by $C(W)\e$. Finally, the second term is easily bounded by $C\e$. Consequently, we have 
\begin{equation*}
    J_\e(v_\e;A_{1-\e,1})\leq C(W)\e.
\end{equation*}
Combining this with the minimality of $u_\e$ implies
\begin{equation}\label{good tr: upper}
\begin{split}
    &J_\e(u_\e,B_1)\leq J_\e(v_\e,B_1)\\
    &\qquad\leq C\e+(|P_1P_3|+|P_2P_4|)(1-\e)+C\e\\
    &\qquad \leq |P_1P_3|+|P_2P_4|+C\e.
\end{split}
\end{equation}

For the lower bound, we extend $u_\e$ to
$B_{1+\e}$ by the interpolation from $g$ to
$h(\f{\cdot}{1+\e})$, i.e. denote by $w_\e$ the function 
\begin{equation}\label{def:w_e}
    w_\e(re^{i\theta}):=\begin{cases}
        u_\e, &\text{ on } B_1,\\
        \f{r-1}{\e} h(\theta)+\f{1+\e-r}{\e} g(\theta),& r\in(1,1+\e),\\
        h(\theta), & \text{ on }\pa B_{1+\e}.
    \end{cases}
\end{equation}
Then we have 
\begin{equation*}
 J_\e(w_\e;A_{1,1+\e}) \leq C\e.
\end{equation*}
Let $\tilde{w}_\e$ minimize the energy in $B_{1+\e}$ with the same trace as $w_\e$.  The lower bound in Corollary \ref{corol: lower bound}, after
scaling by $1+\e$, yields
\begin{equation*}
     J_\e(\tilde{w}_\e;B_{1+\e})
 \ge |P_1P_3|+|P_2P_4|-C\e^{1/3},
\end{equation*}
which directly gives 
\begin{equation}
\begin{split}
 J_\e(u_\e;B_1)
 &\ge J_\e(w_\e;B_{1+\e})-C\e  \\
 &\ge J_\e(\tilde{w}_\e;B_{1+\e})-C\e    \\
 &\ge |P_1P_3|+|P_2P_4|-C\e^{1/3}.   \label{good tr:lower}
 \end{split}
\end{equation}
\eqref{good tr: upper} and \eqref{good tr:lower} together imply \eqref{ene bdd 3}.

It remains to transfer the localization construction in Section \ref{sec:min unit disk}. We denote again by $\{D_i\}_{i=1}^4$ the four partition of $B_1$ determined by $\G$, such that 
\begin{equation*}
    \pa D_i\cap \pa B_1=\arc{P_iP_{i+1}},\quad \g_{i}\subset \pa D_i,\quad \text{for every }i.
\end{equation*}

Invoking the construction in Lemma \ref{lem:4 phase radii} of four line segments $\{\Sigma_i\}_{i=1}^4$, one finds that all the key steps, such as slicing, control of bad sets and horizontal-vertical energy balance,  use the minimality only through the upper energy bound.  Therefore these steps apply to $w_\e$ defined in \eqref{def:w_e},  which is not minimizing, but has exact outer trace $h$ and satisfies the upper bound 
\begin{equation*}
J_\e(w_\e,B_{1+\e})=J_\e(u_\e,B_1)+J_\e(w_\e,A_{1,1+\e})\leq |P_1P_3|+|P_2P_4|+C\e.    
\end{equation*}
Therefore, there are four segments $\{\Sigma_{i}\}_{i=1}^4$ satisfying properties in \eqref{property of Sigma i} and 
\begin{equation*}
    |u_\e-a_i|\leq C\e^{1/12}\quad \text{ on }\Sigma_i, \ \ \text{for every }i.
\end{equation*}

For each \(i\), one can argue in the same way as in the proof of Proposition \ref{prop: loc diff int} to obtain two straight segments $\mathfrak{L}^{i\,i+1}$ and $\mathfrak{L}^{i\,i-1}$ in $D_i$, parallel to the two arms of $\pa D_i\cap B_1$, on which
\begin{equation*}
 | u_\e-a_i|\leq K\bar\delta.
\end{equation*}
The distance of each straight segment to the corresponding arm of $\pa D_i\cap B_1$ is at most $C(W)\e^{1/12}$. $\mathfrak{L}^{i\,i+1}$ and $\mathfrak{L}^{i\,i-1}$ have a common endpoint lying in $B_{C\e^{1/12}}(O_\G)\cap D_i$. 

Let $\widehat\g_i\subset\g_i$ be the boundary subarc between $\mathfrak{L}^{i\,i+1}\cap \pa B_1$ and $\mathfrak{L}^{i\,i-1}\cap \pa B_1$, and let $\widehat D_i\subset D_i\cap B_1$ be the region bounded by $\mathfrak{L}^{i\,i+1}\cup\mathfrak{L}^{i\,i-1}\cup\widehat\gamma_i$. It is easy to see from the properties of good boundary data that
\begin{equation*}
   |u_\e-a_i|\le K\bar\delta \quad\text{ on }\pa\widehat D_i.  
\end{equation*}
Invoking Proposition \ref{prop:maxp} implies 
\begin{equation}\label{good-phase-decay}
 |u_\e-a_i|\le K\bar\delta\quad\text{ in }\widehat D_i,
 \qquad
 |u_\e(z)-a_i|
 \le C\exp\left(-\frac{c\dist(z,\pa\widehat D_i)}\e\right).
\end{equation}
Because $|O_\G|\leq \f1{10}$, the annulus $A_{1/4,3/4}(O_\G)$ stays away from $\pa B_1$.  There is $C_0=C_0(W)$ such that
\begin{equation}\label{shrink phase annulus}
 A_{1/4,3/4}(O_\G)\setminus N_{2C_0\e^{1/12}}(\Gamma)\subset\bigcup_{i=1}^4\left\{z\in\widehat D_i:\dist(z,\pa\widehat D_i)\ge C_0\e^{1/12}\right\}.
\end{equation}
For fixed $\gamma\in(0,K\bar\delta)$,
 when $\e$ is sufficiently small, the exponential decay \eqref{good-phase-decay} and \eqref{shrink phase annulus} guarantee \eqref{annulus loc dif int}, taking $C(\g,W)=2C_0$. This completes the proof of Lemma \ref{lem:fine prop good tr}. 
 
\end{proof}

\begin{lemma}\label{lem:good tr annulus ene}
    Let $u_\e$ minimize $J_\e(\cdot,B_1)$ subject to a good trace $u_\e|_{\pa B_1}=g$ satisfying \eqref{good partition}--\eqref{good center}. And $\G$ and $O_\G$ be the associated cross and center. Then for every pair of radii $(r,\,R)$ such that $\f14\leq r<R\leq \f34$, it holds that 
\begin{equation}\label{ene est: annulus}
    \left|J_\e(A_{r,R}(O_\G))-4(R-r)\right|\leq C\e^{1/12},
\end{equation}
where $C=C(W)$ is independent of $\e$.
\end{lemma}

\begin{proof}
    Fix $(r,R)$ such that $\f14\leq r<R\leq \f34$. \eqref{good-phase-decay} and \eqref{shrink phase annulus} imply that for any $z\in \widehat D_i$ such that $\min\{\dist(z,\G), \,\dist(z,\pa B_1)\}\ge 2C_0\e^{1/12}$, 
    \begin{equation}\label{dist on overhat D}
        |u_\e(z)-a_i|\leq C\exp\left( -C \e^{-\f{11}{12}} \right),
    \end{equation}
    where the value of $C$ may vary, but only depends on $W$. 

    Focus on one arm $O_\G P_i$, set
    \begin{equation*}
        l_i:=|O_\G P_i|,\quad \tau_i:= \f{\overrightarrow{O_\G P_i}}{l_i},\quad \nu_i:= \f{\pi}{2} \text{ rotation of }\tau_i.
    \end{equation*}
    We consider the thin rectangle with length parallel to $O_\G P_i$:
    \begin{equation*}
        \mathfrak{R}_i:= \{O_\G+t\tau_i+s\nu_i: t\in (R,l_i-6C_0\e^{1/12}), s\in (-2C_0\e^{1/12},2C_0\e^{1/12})\}. 
    \end{equation*}
    It is not hard to verify by \eqref{dist on overhat D} that for every $i$,
    \begin{align*}
        &\qquad\qquad \mathfrak{R}_i\subset (B_1\setminus B_R(O_\G))\cap N_{2C_0\e^{1/12}}(\G),\\
        & |u_\e(z)-a_i|\leq C\exp\left( -C \e^{-\f{11}{12}} \right),\quad z\in \{O_\G+t\tau_i+2C_0\e^{1/12}\nu_i: t\in (R,l_i-6C_0\e^{1/12})\},\\
        & |u_\e(z)-a_{i-1}|\leq C\exp\left( -C \e^{-\f{11}{12}} \right),\quad z\in \{O_\G+t\tau_i-2C_0\e^{1/12}\nu_i: t\in (R,l_i-6C_0\e^{1/12})\}.
    \end{align*}
    Thus,
    \begin{equation*}
        J_\e(u_\e,\mathfrak{R}_i)\geq (l_i-R-6C_0\e^{1/12})(1-C \exp(-2C\e^{-\f{11}{12}}))\geq (l_i-R)-C\e^{1/12},
    \end{equation*}
    the last inequality needs $\e$ to be small enough such that $C \exp(-2C\e^{-\f{11}{12}}) $ is smaller than $\e^{1/12}$. Consequently we have 
    \begin{equation*}
        J_\e(u_\e, B_1\setminus B_R(O_\G))\geq \sum_{i=1}^4l_i-4R-C\e^{1/12}.
    \end{equation*}

    Inside the disk $B_r(O_\G)$, for every $s\in (C\e^{1/12},r)$, $\pa B_s(O_\G)$ will intersect with all $\widehat D_i$ and therefore costs energy of four transitions. It is straightforward to obtain
    \begin{equation*}
        J_{\e}(u_\e,B_r(O_\G))\geq 4r- C\e^{1/12}.
    \end{equation*}
    Finally, 
    \begin{equation*}
        J_\e(u_\e,A_{r,R}(O_\G))\leq \sum_{i=1}^4 l_i+C\e- J_\e(u_\e, B_1\setminus B_R(O_\G))- J_{\e}(u_\e,B_r(O_\G))\leq 4(R-r)+C\e^{1/12}. 
    \end{equation*}
    On the other hand, the lower bound 
    \begin{equation*}
        J_\e(u_\e,A_{r,R}(O_\G))\geq 4(R-r)-C\e^{1/12}
    \end{equation*}
    follows immediately from the fact that every $\pa B_s(O_\G)$ passes through all $\widehat D_i$ for $s\in(r,R)$. The proof is completed. 
\end{proof}
\begin{figure}[htb!]
\centering
\begin{tikzpicture}[scale=2]
\definecolor{primary}{HTML}{08457E}
\definecolor{secondary}{HTML}{F209B1}
\definecolor{interfacegray}{HTML}{68737D}
  \def\ang{42}
  \def\rad{1.42}
  \def\rsel{0.88}
  \coordinate (O) at (0,0);
  \coordinate (P1) at ({\rad*cos(-\ang)},{\rad*sin(-\ang)});
  \coordinate (P2) at ({\rad*cos(\ang)},{\rad*sin(\ang)});
  \coordinate (P3) at ({\rad*cos(180-\ang)},{\rad*sin(180-\ang)});
  \coordinate (P4) at ({\rad*cos(180+\ang)},{\rad*sin(180+\ang)});
  \fill[primary!6] (O) circle (1.1);
  \fill[white] (O) circle (0.7);
  \draw[primary!30] (O) circle (0.7);
  \draw[primary!30] (O) circle (1.1);
  \begin{scope}
    \clip (O) circle (\rad);
    \draw[primary!12,line width=12pt] (P1)--(P3);
    \draw[primary!12,line width=12pt] (P2)--(P4);
    \draw[primary,densely dashed] (P1)--(P3);
    \draw[primary,densely dashed] (P2)--(P4);
  \end{scope}
  \draw[interfacegray] (O) circle (\rad);
  \draw[secondary,very thick] (O) circle (\rsel);
  \foreach \a in {-42,42,138,222}
    {\fill[primary] ({\rsel*cos(\a)},{\rsel*sin(\a)}) circle (1.4pt);}
  \coordinate (Q1) at ({\rsel*cos(-38)},{\rsel*sin(-38)});
  \coordinate (Q2) at ({\rsel*cos(47)},{\rsel*sin(47)});
  \coordinate (Q3) at ({\rsel*cos(145)},{\rsel*sin(145)});
  \coordinate (Q4) at ({\rsel*cos(224)},{\rsel*sin(224)});
  \foreach \a in {-38,47,145,224}
    {\draw[secondary,line width=4pt]
      ({\rsel*cos(\a-2.8)},{\rsel*sin(\a-2.8)})
      arc[start angle={\a-2.8},end angle={\a+2.8},radius=\rsel];}
  \draw[secondary,thick] (Q1)--(Q3) (Q2)--(Q4);
  \fill[secondary] (Q1) circle (1.8pt);
  \fill[secondary] (Q2) circle (1.8pt);
  \fill[secondary] (Q3) circle (1.8pt);
  \fill[secondary] (Q4) circle (1.8pt);
  \node[primary] at (-0.88,1.2) { $O(\e^{1/12})$ neighborhood};
  \node[primary] at (-0.98,-0.73) { $\Gamma$};
  \node[secondary] at (0.6,-0.32) { $\widetilde\Gamma$};
  \node[secondary] at (0.96,0.12) { $\partial B_\rho(O_\Gamma)$};
\end{tikzpicture}
\caption{Good radius gives an approximating cross $\tilde{\G}$ at a smaller scale}
\label{good radius}
\end{figure}
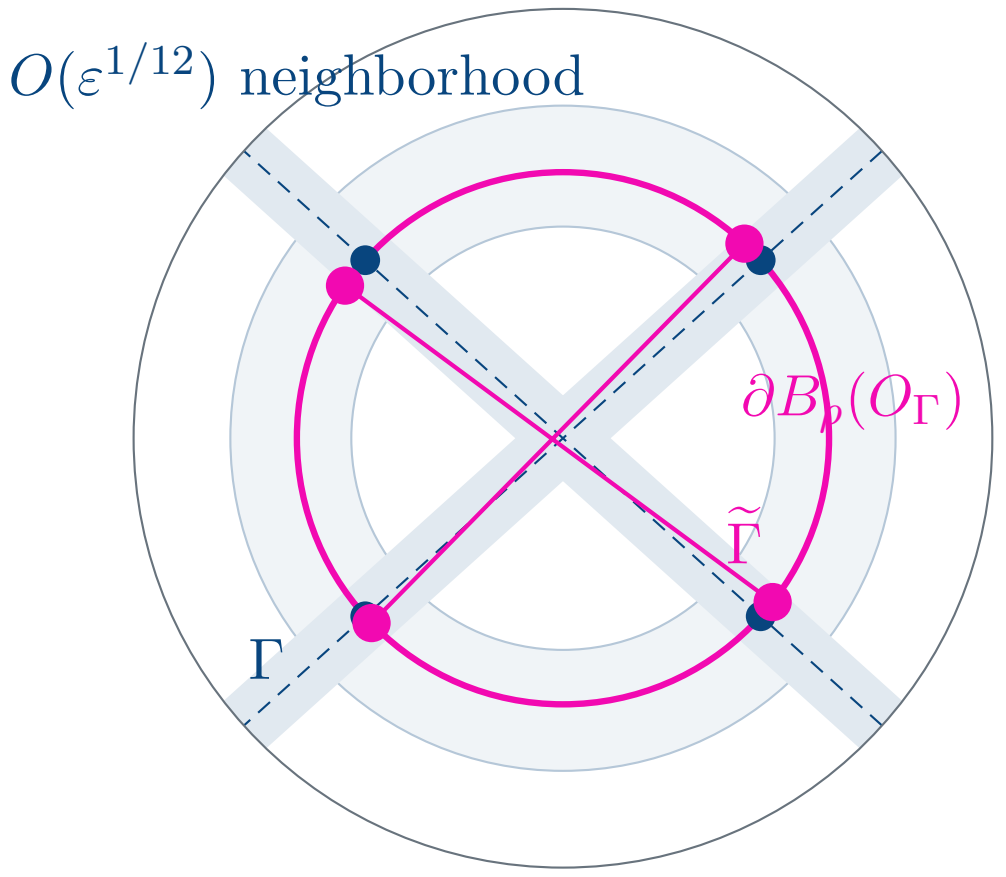

\begin{proposition}\label{prop: good radius}
There exist constants $\e_0>0$, $C>0$ depending only on $W$ with the following property. 

Take $\e<\e_0$ and let $u_\e$ be a minimizer of $J_\e(\cdot, B_1)$ as in Lemma \ref{lem:fine prop good tr} with $u_\e|_{\pa B_1}$ being a good trace satisfying \eqref{good partition}--\eqref{good center}. We further assume the following stronger version of \eqref{good cross angle} to leave some margin for minimality:
   \begin{equation}\label{angle margin}
   \angle P_iO_\G P_{i+1}\in  \bigg(2\arccos(\f{s}{2})+\mu_0, \ 2\arcsin(\f{s}{2})-\mu_0 \bigg),   
   \end{equation}
   for a constant $\mu_0\in\left(0, \f{\arcsin(\f{s}{2})-\arccos(\f{s}{2})}{2}\right)$.

Then there is $\rho\in[1/3,1/2]$ such that $u_\e|_{\pa B_\rho(O_\Gamma)}$, after the rescaling to unit disk:
\begin{equation*}
 \zeta=\frac{z-O_\Gamma}{\rho} \qquad v_{\tilde\e}(\zeta)=u_\e(O_\Gamma+\rho\zeta),\qquad
 \tilde\e=\f\e\rho,
\end{equation*}
is again a good trace on $\pa B_1$ satisfying \eqref{good partition}--\eqref{good center}. 
   
$v_{\tilde{\e}}$ has the new associated small boundary gaps $\{\tilde{G}_i\}_{i=1}^4$ and phase arcs $\{\tilde{\g}_{i}\}_{i=1}^4$, as well as the new approximating cross $\tilde{\G}$ determined by midpoints of $\tilde{G}_i$. It holds that 
   \begin{equation}\label{est:angle change 1 step}
       d_{angle}(\G,\tilde{\G}):= \max\{\angle( P_1P_3, \tilde{P}_1\tilde{P}_3), \ \ \angle(P_2P_4, \tilde{P}_2\tilde{P}_4) \}\leq C\e^{1/12},
   \end{equation}
   \begin{equation}\label{est:center change 1 step}
       |O_{\tilde{\G}}|\leq C\e^{1/12},
   \end{equation}
   for some constant $C=C(W)$.
\end{proposition}

\begin{proof}
For an arc $S$ on a circle centered at $O_\G$, write
\begin{equation*}
\mathcal E_\e(S):=\int_S\left(\frac\e2|\pa_\tau u_\e|^2+\frac1\e W(u_\e)\right)\,d\ch.
\end{equation*}
By Lemma \ref{lem:good tr annulus ene} and Fubini, there is $\rho\in[1/3,1/2]$ such that, for sufficiently small $\e$,
\begin{equation}\label{circle-energy}
 \mathcal E_\e\left(\pa B_\rho(O_\G)\right)
 \leq 4+C(W)\e^{1/12}\leq 4+C\bar\delta^2,
\end{equation}
Here, the constant $C$ before $\bar\delta^2$ is the same as in \eqref{est:endpoint close to ai} and \eqref{good energy}.

\eqref{shrink phase annulus} and \eqref{good-phase-decay} give four disjoint arcs
\begin{equation*}
T_i:=\pa B_\rho(O_\G)\cap N_{C\e^{1/12}}(O_\G P_i), \quad \ch(T_i)\leq C\e^{1/12},\qquad i=1,...,4,    
\end{equation*}
such that
\begin{equation*}
 |u_\e-a_i|\leq\f{\bar\delta}{2}\quad\text{on the arc between }T_i\text{ and }T_{i+1}.
\end{equation*}

Parameterizing $\mathbb{S}^1$ by $\theta\in [0,2\pi)$, we assume 
\begin{equation*}
    T_i=  (\theta_i^1,\theta_i^2),\qquad  \theta_i^1<\theta_i^2<\theta_i^1+C\e^{1/12}.
\end{equation*}
Then by definition of $T_i$ we have 
\begin{equation*}
    |u_\e(\theta_i^1)-a_{i-1}|\leq \f{\bar\delta}{2},\qquad |u_\e(\theta_i^2)-a_{i}|\leq \f{\bar\delta}{2}
\end{equation*}
Define
\begin{align*}
 \theta_i^-&:=\sup\{\theta\in (\theta_i^1,\theta_i^2):|u_\e(\theta)-a_{i-1}|\le\bar\delta\},\\
 \theta_i^+&:=\inf\{\theta\in(\theta_i^-, \theta_i^2):|u_\e(\theta)-a_i|\le\bar\delta\}, \qquad G_i^\rho:=[\theta_i^-,\theta_i^+].
\end{align*}
Continuity implies
\begin{equation*}
 |u_\e(\theta_i^-)-a_{i-1}|=|u_\e(\theta_i^+)-a_i|=\bar\delta, \qquad G_i^\rho\subset\overline{T_i},  
\end{equation*}
which further yields
\begin{equation}\label{eq:gap-energy}
 \mathcal E_\e(G_i^\rho)\ge1-C\bar\delta^2,
\end{equation}
where $C$ is again the same constant as in \eqref{good energy}.

We claim that
\begin{equation}\label{G i rho measure}
    \min_{j=1,2,3,4} |u_\e-a_j|\geq \bar{\delta}\quad \text{ on }G_{i}^\rho.
\end{equation}
Indeed, $u_\e$ cannot be $\bar\delta$-close to $a_i$ or $a_{i-1}$ by the definition of $G_i^\rho$. For other $j\neq i,i-1$, if $u_\e$ visits $\bar\delta$ neighborhood of $a_j$ on $G_i^\rho$, there will be an extra phase transition created on $\pa B_\rho(O_\G)$, which together with the energy contribution from $G_k^\rho$ for $k\neq i$ yields a contradiction with \eqref{circle-energy}. It follows directly from \eqref{G i rho measure} that
\begin{equation}\label{eq:gap-length}
 \ch(G_i^\rho)\leq \f{\e}{c\bar\delta^2}\mathcal E_\e(G_i^\rho) \leq \left(\f{5}{c\bar\delta^2}\right)\e=C_g\e.
\end{equation}

Let $\gamma_i^\rho$ be the complementary arc between $G_i^\rho$ and $G_{i+1}^\rho$.  Subtracting the four transitional costs on $\{G_i^\rho\}_{i=1}^4$ from \eqref{circle-energy} gives
\begin{equation}\label{residual}
 \sum_{i=1}^4\mathcal E_\e(\g_i^\rho)\leq 5C\bar\delta^2.
\end{equation}
If $|u_\e-a_i|>K\bar\delta$ somewhere on $\g_i^\rho$, then on $\g_i^\rho$, $u_\e$ needs to travel from $\bar\delta$ neighborhood of $a_i$ to the outside of $K\bar\delta$ neighborhood of $a_i$, which costs energy at least
\begin{equation*}
    c_W(K\bar\delta)^2-C_W\bar\delta^2\geq 8C\bar\delta^2,
\end{equation*}
contradicting \eqref{residual}.  Therefore,
\begin{equation}\label{pure arcs}
 |u_\e-a_i|\le K\bar\delta\quad\text{on }\gamma_i^\rho,\quad \text{for every }i.
\end{equation}

Now set
\begin{equation*}
   \tilde G_i=\f{G_i^\rho-O_\G}{\rho},\qquad\tilde\g_i=\frac{\g_i^\rho-O_\G}{\rho}, \qquad v_{\tilde\e}(\zeta)=u_\e(O_\G+\rho\zeta), \qquad \tilde\e=\f\e\rho. 
\end{equation*}

The estimates \eqref{G i rho measure} and \eqref{pure arcs} are unchanged by this normalization, while
\begin{align*}
 &\qquad \ch(\tilde G_i)\le C_g\frac\e\rho=C_g\tilde\e,\\
\int_{\pa B_1} &\left(\frac{\tilde\e}{2}|\pa_\tau v_{\tilde\e}|^2+\frac1{\tilde\e}W(v_{\tilde\e})\right)\,d\ch=\mathcal E_\e\left(\pa B_\rho(O_\G)\right)
 \leq 4+C\bar\delta^2.
\end{align*}
Thus \eqref{good partition}--\eqref{good energy} hold for the new trace.

It remains to estimate the angle.  Let $\tilde P_i$ be the midpoint of $\tilde G_i$, and set
\begin{equation*}
 \tau_i:=\frac{P_i-O_\Gamma}{|P_i-O_\Gamma|}\in \pa B_1.
\end{equation*}
By the choice of $G_i^\rho$ we have 
\begin{equation}\label{boundary pts close}
 \max_i \dist(\tilde P_i,\tau_i)\leq C\e^{1/12}.
\end{equation}
Since $O_\G\in P_1P_3\cap P_2P_4$, $\tau_1=-\tau_3$, $\tau_2=-\tau_4$. \eqref{boundary pts close} and elementary geometry implies \eqref{est:angle change 1 step} and \eqref{est:center change 1 step}. Finally, \eqref{est:angle change 1 step} and \eqref{angle margin} together indicate that $\tilde{\G}$ satisfies the condition \eqref{good cross angle} as long as we choose $C\e^{1/12}\ll \mu_0$. Thus $v_{\tilde{\e}}|_{\pa B_1}$ is again a good trace. The proof is complete.  
\end{proof}

\section{Proof of Theorem \ref{thm:main}}\label{sec:min thm proof}

Let $\G$ be the prescribed cross, with $O_\G=0$, and let $\{D_i\}_{i=1}^4$ be the four-partition of $\BR^2$ determined by $\G$. $\{P_i\}_{i=1}^4$ are four boundary points such that $ D_i\cap \pa B_1=\arc{P_iP_{i+1}}$. We write
\begin{equation*}
 u_\G=\sum_{i=1}^4a_i\mathbf1_{D_i}.
\end{equation*}
Choose $\mu_0>0$ sufficiently small that
\begin{equation}\label{iteration: ini margin}
 \angle P_iO_\Gamma P_{i+1}
 \in\left(2\arccos\frac s2+4\mu_0,\,
          2\arcsin\frac s2-4\mu_0\right),
 \qquad i=1,...,4.
\end{equation}
Throughout this section, constants $c,C>0$ depend only on $W$, while the fixed small threshold $\e_*>0$ chosen below may also depend on $\mu_0$.

Assume $u_\e$ minimizes $J_\e(\cdot,B_1)$ with the trace $g_{\e,\G}$ defined in Section \ref{sec:min unit disk} satisfying \eqref{est:trace}.  By Lemma \ref{lem:upper bdd} and Lemma \ref{lem:good tr annulus ene} we have 
\begin{equation}\label{ini upper bdd}
 J_\e(u_\e,B_1)\leq4+C\e,
 \qquad \|u_\e\|_{L^\infty(B_1)}\leq C.
\end{equation}
\begin{equation}\label{ini annular bound}
 \left|J_\e(u_\e,A_{r,R})-4(R-r)\right|
 \leq C\e^{1/12},\qquad \f14\leq r<R\leq \f34.
\end{equation}

We carry out the following iteration in the original coordinates.  At step $j$, let $B_{r_j}(O_j)$ be the current disk and let $\G_j$ be its associated cross centered at $O_{j+1}$. Initially,
\begin{equation*}
    r_0=1,\qquad O_0=O_1=0,\qquad \G_0=\G. 
\end{equation*}
At each step $j$, we also define the rescaled map and its cross:
\begin{equation}\label{iter rescale}
\begin{split}
 &\qquad v_j(\zeta):=u_\e(O_j+r_j\zeta),\qquad \tilde{\G}_j:=\f{\G_j-O_j}{r_j},\\
 &\quad J_\e(u_\e,O_j+r_j\Omega)
 =r_jJ_{\e/r_j}(v_j,\Omega),\qquad \e_j:=\f{\e}{r_j}.
 \end{split}
\end{equation}
Suppose $v_j$ is minimizing on $B_1$, and $v_j|_{\pa B_1}$ is good with parameter $\e_j$. More precisely, we require that $v_j|_{\pa B_1}$ satisfies \eqref{good partition}--\eqref{good energy} and \eqref{angle margin}.  Apply Proposition \ref{prop: good radius} at this scale to obtain
\begin{equation}\label{disk renewal}
 r_{j+1}=\rho_jr_j\quad
 \text{ for some }\rho_j\in[\f13,\f12],\qquad O_{j+1}:=O_{\G_j}, 
\end{equation}
such that $u_\e|_{\pa B_{r_{j+1}}(O_{j+1})}$, after rescaling back to the unit disk, is a good trace again.  And the renewed cross $\G_{j+1}\subset B_{r_{j+1}}(O_{j+1})$ satisfies
\begin{align}
 \label{it 1 angle diff}&d_{angle}(\G_{j+1},\G_j)
 \leq C\e_j^{1/12}= C\left(\f\e{r_j}\right)^{1/12},\\
 \label{it 1 center diff}&|O_{\Gamma_{j+1}}-O_{\Gamma_j}|
 \leq Cr_j\e_j^{1/12}= Cr_j\left(\f\e{r_j}\right)^{1/12}
 \leq C\e^{1/12}r_j^{11/12}.     
\end{align}

The geometric decay of $r_j$ implies
\begin{equation}\label{forward sum}
 \sum_{k=0}^{j-1}\left(\f\e{r_k}\right)^{1/12}\leq\left(\frac\e{r_j}\right)^{1/12}\left(\sum_{k=1}^{\infty}2^{-k/12}\right)\leq C\left(\frac\e{r_j}\right)^{1/12}.
\end{equation}
Consequently,
\begin{equation}\label{cross sum}
 d_{angle}(\Gamma_j,\Gamma)
 \le C\left(\frac\e{r_j}\right)^{1/12}.
\end{equation}
Choose $\e_*>0$ sufficiently small such that the total change of the opening angle for the cross $\G_j$ compared to the initial cross $\G$ is less than $\mu_0$ whenever $\e_j\leq\e_*$, and that all smallness requirements on $\e$ in Lemma \ref{lem:fine prop good tr} and Proposition \ref{prop: good radius} hold for $\e_j<\e_*$.  Thus \eqref{iteration: ini margin} would ensure the angle requirement \eqref{angle margin} for next iteration at every such step. 

We stop the iteration at the first index $N=N(\e)$ for which
\begin{equation}\label{stopping step}
 \f\e{r_N}\geq\f{\e_*}{3}.
\end{equation}
By minimality of $N$, we have
\begin{equation*}
 \frac\e{\e_*}<r_N\le\frac{3\e}{\e_*}.
\end{equation*}

Let
\begin{equation*}
 z_\e:=O_{\G_N}.
\end{equation*}
Summing \eqref{it 1 center diff} from step $j$ to $N$ yields
\begin{equation}
    |z_\e-O_{\G_j}|\leq C\e^{1/12}\sum_{k=j}^{N-1}r_k^{11/12} \leq C\e^{1/12}r_j^{11/12}.                           \label{tail center diff} 
\end{equation}
In particular, taking $j=0$ implies 
\begin{equation}\label{terminal center}
 |z_\e|\le C\e^{1/12}.
\end{equation}

We next estimate how $u_\e$ resembles $u_\G$ after recentering at $z_\e$. Denote by 
$$
u_{\Gamma_j}:=\sum_i a_i\mathbf1_{D_i^j}
$$
the sharp four-phase map  associated with the $j$-step cross $\G_j$, here $\{D_i^j\}_{i=1}^4$ is the four-partition determined by $\G_j$.  For every $j\geq 1$, \eqref{good-phase-decay} and \eqref{shrink phase annulus}, after scaling, imply
\begin{equation} \label{annulus L1}
\begin{split}
 &\int_{A_{r_j/4,\,3r_j/4}(O_{\G_j})} |u_\e-u_{\Gamma_j}|\,dz \\
 &\qquad\leq C r_j^2\left(\frac\e{r_j}\right)^{1/12} +C r_j^2\exp\left[-c\left(\frac{r_j}{\e}\right)^{11/12}\right] \leq C\e^{1/12}r_j^{23/12}.   
\end{split}
\end{equation}
Here the first term is the area of the $O(r_j\e_j^{1/12})$ neighborhood around $\G_j$, and the second term records the exponential decay to corresponding phases outside the diffuse interface. For $j=0$, the same estimate follows from \eqref{est: exponential}.

By \eqref{tail center diff} and \eqref{cross sum} and elementary geometric estimates on the translation of disks or annuli by a small distance, we have 
\begin{equation}
 \left|A_{r_j/4,\,3r_j/4}(z_\e)\,\triangle\, A_{r_j/4,\,3r_j/4}(O_{\Gamma_j})\right| \leq Cr_j|z_\e-O_{\Gamma_j}|\leq C\e^{1/12}r_j^{23/12}.                         \label{annulus translation}
\end{equation}

\begin{equation} \label{comp Gj, G}
\begin{split}
 &\int_{B_{r_j}(z_\e)}|u_{\G_j}(z)-u_\G(z-z_\e)|\,dz \\
 \leq & \int_{B_{r_j}(z_\e)}\bigg(|u_{\G_j}(z)-u_\G(z-O_{\G_j})| + |u_\G(z-O_{\G_j})-u_\G(z-z_\e)|\bigg)\,dz\\
 \leq & \ C r_j^2d_{angle}(\Gamma_j,\Gamma)+ C r_j|O_{\Gamma_j}-z_\e|\leq C\e^{1/12}r_j^{23/12}.                        
\end{split}
\end{equation}
Combining \eqref{annulus L1}--\eqref{comp Gj, G}, we obtain
\begin{equation}\label{comp annulus ue uG}
 \int_{A_{r_j/4,\,3r_j/4}(z_\e)} |u_\e(z)-u_\G(z-z_\e)|\,d z\leq C\e^{1/12}r_j^{23/12},\qquad 0\leq j\leq N.
\end{equation}
Given an arbitrary $r\in[r_N,\f12]$, choose the largest integer $k$ such that $\f34r_k\geq r$. \eqref{comp annulus ue uG} further yields
\begin{equation}\label{Br L1}
\begin{split}
 &\int_{B_r(z_\e)}|u_\e(z)-u_\Gamma(z-z_\e)|\,dz\\
 &\leq \sum_{j=k}^N \int_{A_{r_j/4,\,3r_j/4}(z_\e)} |u_\e(z)-u_\G(z-z_\e)|\,dz +\int_{B_{\f{r_N}{4}}(z_\e)} |u_\e(z)-u_\Gamma(z-z_\e)|\,dz \\
 &\leq C\e^{1/12}\sum_{j=k}^N r_j^{23/12}+Cr_N^2 \leq C\e^{1/12}r^{23/12}+C\left(\f\e{\e_*}\right)^2.      
\end{split}
\end{equation}

Define the blow-up maps
\begin{equation}\label{eq:blowup-maps}
 U_\e(x):=u_\e(z_\e+\e x),\qquad x\in\f{B_1-z_\e}{\e}.
\end{equation}
Then 
\begin{equation*}
    \Delta U_\e=W_u(U_\e),\qquad\|U_\e\|_{L^\infty}\leq C.
\end{equation*}
It follows from interior elliptic estimates that, along a subsequence,
\begin{equation}\label{entire-limit}
U_\e\rightarrow U \quad\text{ in }C^2_{loc}(\BR^2;\BR^2),\qquad \Delta U=W_u(U).
\end{equation}
It is straightforward to verify that $U$ inherits the minimality from $u_\e$.

Finally we identify the blow-down limit \eqref{eq:blowdown} of $U$ at infinity. It suffices to prove
\begin{equation*}
    \lim\limits_{R\to \infty}\int_{B_1} \left|U(Rz)-u_\G(z)\right|\,dz=0.
\end{equation*}
Take $r=\e R$ in \eqref{Br L1}. After blow-up we have  
\begin{equation*}
\int_{B_R}|U_\e(z)-u_\Gamma(z)|\,dz\leq C R^{23/12}+C\e_*^{-2}.
\end{equation*}

Passing to the limit $\e\to 0$ yields
\begin{equation}\label{eq:large-ball-L1}
 \int_{B_R}|U-u_\Gamma|\,dz\leq C R^{23/12}+C\e_*^{-2}, \qquad R\geq4/\e_*,
\end{equation}
which is equivalent to
\begin{equation*}
    \int_{B_1}|U(Rz)-u_\Gamma(z)|\,dz \leq C R^{-1/12}+C\e_*^{-2}R^{-2}\to 0  \qquad\text{as }R\to\infty.
\end{equation*}

This completes the proof of Theorem \ref{thm:main}.

\bibliographystyle{acm}
\bibliography{cross.bib}

\end{document}